\documentstyle{article}

\newcommand{\N}{{\rm I\kern-.5ex N}}
\newcommand{\Z}{{\sf \vrule height 1.55ex depth-1.2ex
width.03em\kern-.11em Z
\kern-.9ex Z\kern-.11em\vrule height 0.3ex depth0ex width.03em}}
\newcommand{\Q}{{\rm\kern.2ex\vrule height1.55ex depth-.05ex
width.03em\kern-.7ex Q}}
\newcommand{\R}{{\rm I\kern-.5ex R}}
\newcommand{\Rvar}{{\rm I\kern-.5ex R}}
\newcommand{\C}{{\rm\kern.3ex\vrule height1.55ex depth-.05ex
width.03em\kern-.7ex C}}

\newcommand{\hoed}{\hat{\rule{0ex}{1.7ex}}\,}
\newcommand{\spat}{\hspace{4ex}}

\newcommand{\ahh}{\hat{A}  
\hspace{-.55ex}\hat{\rule{0ex}{2.0ex}}\hspace{.55ex}}
\newcommand{\dhh}{\hat{\de}  
\hspace{-.95ex}\hat{\rule{0ex}{2.05ex}}\hspace{.95ex}}

\newcommand{\flip}{ \chi }

\newcommand{\cst}{C$^*$}

\newcommand{\ah}{\hat{A}}
\newcommand{\deh}{\hat{\Delta}}
\newcommand{\psih}{\hat{\psi}}
\newcommand{\vfih}{\hat{\vfi}}

\newcommand{\od}{\odot}
\newcommand{\ot}{\otimes}

\newcommand{\om}{\omega}
\newcommand{\io}{\iota}
\newcommand{\vfi}{\varphi}
\newcommand{\vep}{\varepsilon}

\newcommand{\sde}{\delta}
\newcommand{\de}{\Delta}

\newcommand{\th}{\theta}

\newcommand{\text}[1]{\mbox{#1}}

\newcommand{\qed}{\ \hfill \rule{2mm}{2mm}}
\newenvironment{demo}{\medskip\noindent\bf Proof :\ \
\rm}{\qed\bigskip\par }

\newtheorem{definition}{Definition}[section]
\newtheorem{proposition}[definition]{Proposition}
\newtheorem{lemma}[definition]{Lemma}
\newtheorem{corollary}[definition]{Corollary}
\newtheorem{remark}[definition]{Remark}
\newtheorem{theorem}[definition]{Theorem}

\newtheorem{notation}[definition]{Notation}
\newtheorem{result}[definition]{Result}

\begin{document}

\begin{center}
\LARGE\bf
Examining the dual of an algebraic quantum group.

\end{center}

\bigskip

\begin{center}
\rm J. Kustermans  \footnote{Research Assistant of the
National Fund for Scientific Research (Belgium)}

Institut for Matematik og Datalogi

Odense Universitet

Campusvej 55

5230 Odense M

Denmark

\bigskip

April 1997
\end{center}

\subsection*{Abstract}

In the first part of this paper, we implement the multiplier algebra
of the dual of an algebraic quantum group $(A,\de)$ (see \cite{VD1})
as a subset of the space of linear functionals on $A$.

In a second part, we construct the universal corepresentation and use
it to prove a bijective corespondence between corepresentations of
$(A,\de)$ and homomorphisms on its dual.

\section*{Introduction}

In \cite{VD1}, A. Van Daele introduces and investigates a class of
algebraic quantum groups. These algebraic quantum groups are, loosely
speaking, non-unital Hopf algebras which have a left invariant
functional.

These algebraic quantum groups behave very well, as can be found in
\cite{VD1}. The left invariant functional is faithful and unique. It
also satifies some weak KMS property. You can also introduce the
modular function of  such an algebraic quantum group.

It is also possible to construct the dual of an algebraic quantum
group and get again an algebraic quantum group. An overview of the
most important properties can be found in the first section.

\medskip

The canonical examples of algebaic quantum groups  are the compact and
algebraic quantum groups. However, the double construction of
Drinfel'd guarantees the existence of non-discrete non-compact
algebraic quantum groups.

\medskip

If an algebraic quantum group has a compatible $^*$-structure and a
positive left Haar functional, it is possible to construct
\cst-algebraic quantum groups out of them in the sense of Masuda, Nakagami
\& Woronowicz
(see \cite{Kus} and \cite{Kus2}).

\bigskip

Consider an algebraic quantum group $(A,\de)$. Then we get a dual
algebraic quantum group $(\ah,\deh)$ where $\ah$ is a certain subspace
of $A'$. We prove in section 2 that the multiplier algebra $M(\ah)$
can also be considered as a certain subspace of $A'$. This is then
used in section 3 to get a familiar implementation of the
comultiplication $\deh$.

\medskip

In section 5, we prove some useful properties about coreperesentations
on algebraic quantum groups. In section 6, we construct the universal
corepresentation of $(A,\de)$. We get a bijective correspondence
between non-degenerate corepresentations of $(A,\de)$ and
non-degenerate homomorphisms on $\ah$ where this universal
corepresentation serves as the linking bridge.

\bigskip

All vector spaces in this paper are considered over the complex
numbers. For every vector space $V$, we will denote the algebraic dual
by $V'$ and the set of linear operators  on $V$ by $L(V)$.

\medskip

If $V$,$W$ are two vector spaces, then $\flip$ will denote the flip
map $\flip : V \od W \rightarrow W \od V$.

\section{Algebraic quantum groups}

In this  section, we will introduce the notion of an algebraic quantum group
as can be found in \cite{VD1}. Moreover, we will give an overview of the
properties of this algebraic quantum group. The proofs of these results
can be found in the same paper \cite{VD1}. After this section, we will prove 
some further properties about these algebraic quantum groups. We will first
introduce some terminology.

\bigskip

We call an algebra $A$ non-degenerate if and only if we have for every 
$a\in A$ that:
$$(\forall b \in A: a b =0) \Rightarrow a=0 \hspace{1.5cm} \text{ and }
\hspace{1.5cm} (\forall b \in A: b a =0) \Rightarrow a=0 .$$

For a non-degenerate algebra $A$, you can define the multiplier
algebra $M(A)$. This is a unital algebra in which $A$ sits as a
two-sided ideal.

\medskip

If you have two non-degenerate algebras $A,B$ and a multiplicative linear 
mapping $\pi$ from $A$ to $M(B)$, we call $\pi$ non-degenerate if and only
if the vectorspaces $\pi(A) B$ and $B \, \pi(A)$ are equal to $B$. Such a
non-degenerate multiplicative linear map has a unique multiplicative
linear extension to $M(A)$, this extension will be denoted $\overline{\pi}$. 
For every $a \in M(A)$, we define $\pi(a) = \overline{\pi}(a)$.

We have of course similar definitions and results for
anti-multiplicative mappings.

\medskip

For a linear functional $\om$ on a non-degenerate algebra $A$ and any 
$a \in M(A)$ we define the linear functionals $\om a$ and $a \om$ on $A$
such that $(a \om)(x) = \om(x a)$ and $(\om a)(x) = \om(a x)$ for every $x
\in  A$.

\medskip
 
You can find some more information about non-degenerate algebras in the
appendix of \cite{VD4}.

\bigskip

Let $\om$ be a linear functional on an algebra $A$. Then $\om$ is said to
be faithful if and only if we have for  every $a \in A$ that
$$(\forall b \in A: \om(a b)=0) \Rightarrow a=0 \hspace{1.5cm} 
\text{ and} \hspace{1.5cm} (\forall b \in A: \om(b a) =0) \Rightarrow a=0 .$$

\bigskip

We have now gathered the necessary information to understand the
following definition

\begin{definition}
Consider a  non-degenerate algebra $A$ and a non-degenerate homomorphism
$\de$ from $A$ into $M(A \od A)$ such that
\begin{enumerate}
\item $(\de \od \io)\de = (\io \od \de)\de$.
\item The linear mappings $T_1$, $T_2$, $T_3$, $T_4$ from $A \od A$ into
$M(A \od A)$ such that
$$T_1(a \ot b) = \de(a)(b \ot 1) \hspace{1.5cm} \text{ and }
\hspace{1.5cm} T_2(a \ot b) = \de(a)(1 \ot b)$$
$$T_3(a \ot b) = (b \ot 1)\de(a) \hspace{1.5cm} \text{ and }
\hspace{1.5cm} T_4(a \ot b) = (1 \ot b)\de(a)$$
for all $a,b \in A$, are bijections from $A \od A$ to $A \od A$.
\end{enumerate}
Then we call $(A,\de)$ a regular Multiplier Hopf algebra.
\end{definition}

\medskip

\begin{remark} \rm
\begin{itemize}
\item Let $a \in A$ and $\om \in A'$, then $(\om \od \io)\de(a)$ will be
by definition the element in $M(A)$ such that
$$ [(\om \od \io)\de(a)] \, b = (\om \od \io)(\de(a)(1 \ot b))
\hspace{1.5cm} \text{ and }  \hspace{1.5cm}  b \, [(\om \od \io)\de(a)] =
(\om \od \io)((1 \ot b) \de(a)) $$
for every $b \in B$ (if $\de(a)$ would belong to $A \od A$, this
definition of $(\om \od \io)\de(a)$ would be equal to the usual notion
of $(\om \od \io)\de(a)$ \  ).

\item It is not difficult to check that  we have for $\om \in A'$ and
$a,b \in A$ that $(\om b \od \io)\de(a)$ belongs to $A$ and
$(\om b \od \io)\de(a) = (\om \od \io)( (b \ot 1) \de(a) )$.
A similar remark applies for $b \, \om$.

\item Of course, we will use a similar notation for $(\io \od \om)\de(a)$
if $\om \in A'$ and $a \in A$
\end{itemize}
\end{remark}

\medskip

In \cite{VD4}, A.\ Van Daele proves the existence of a unique non-zero
homomorphism $\vep$ from $A$ to $\C$ such that
$$(\vep \od \io)\de = (\io \od \vep)\de =\io \ .$$
He proves moreover the existence of a
unique anti-automorphism $S$ on $A$ such that
$$m(S \od \io)(\de(a)(1 \ot b)) = \vep(a) b \hspace{1.5cm}
{ and } \hspace{1.5cm} m(\io \od S)((b \ot 1)\de(a)) = \vep(a) b $$
for every $a,b \in A$  (here, $m$ denotes the multiplication map from $A
\od A$ to $A$).

As usual, $\vep$ is called the counit and $S$ the antipode of
$(A,\de)$. Furthermore, $\flip(S \od S)\de = \de S$.

\bigskip

Let $\om$ be a linear functional on $A$. We call $\om$ left invariant
(with respect to $(A,\de)$), if and only if $(\io \od \om)\de(a) = \om(a)
1$ for every $a \in A$.
Right invariance is defined in a similar way.

\begin{definition}
Consider a regular Multiplier Hopf algebra $(A,\de)$ such that there
exists a non-zero linear functional $\vfi$ on $A$ which is left
invariant.  Then we call $(A,\de)$ an algebraic quantum group.
\end{definition}

Such a non-zero left invariant linear functional $\vfi$ will be called a
left Haar functional on $(A,\de)$.

\bigskip\bigskip

For the rest of this section, we will fix an algebraic quantum group with
a left Haar functional $\vfi$ on it. We define $\psi$ as the linear
functional $\vfi S$. It is clear that $\psi$ is a non-zero right
invariant linear functional on $A$.

\medskip

An important feature of such an algebraic quantum group is the
faithfulness and uniqueness of left invariant functionals:
\begin{enumerate}
\item Consider a left invariant linear functional $\om$ on $A$, then
there exists a unique element $c \in \C$ such that $\om = c \, \vfi$.
\item Consider a non-zero left invariant linear functional $\om$ on $A$,
then $\om$ is faithful.
\end{enumerate}
In particular, $\vfi$ is faithful.

\medskip

We have of course similar faithfulness and uniqueness results about
right invariant linear functionals.

\bigskip

A first application of this uniqueness result concerns the antipode:
Because $\vfi S^2$ is left invariant, there exists a unique complex
number $\mu$ such that $\vfi S^2 = \mu \vfi$ (in \cite{VD1}, our $\mu$
is denoted by $\tau$!).

\bigskip

In this paper, we will need the following formula:
\begin{equation}
(\io \od \vfi)((1 \ot a)\de(b))
= S(\,(\io \od \vfi)(\de(a)(1 \ot b))\,)  \label{eq1.1} 
\end{equation}
for all $a,b \in A$. A proof of this result can be found in proposition
3.11 of \cite{VD1}.

Using the formula $\flip (S \od S) \de = \de S$, we get the following
form of the formula above.
\begin{equation}
(\psi \od \io)((a \ot 1 \de(b))
= S^{-1}(\,(\psi \od \io)(\de(a)(b \ot 1))\,)  \label{eq1.2}
\end{equation}

\bigskip

Another non-trivial property about $\vfi$ is the existence of a unique
automorphism $\rho$ on $A$ such that $\vfi(a b) = \vfi(b \rho(a))$ for
every $a,b \in A$. We call this the weak KMS-property of $\vfi$ (In
\cite{VD1}, our mapping $\rho$ is denoted by $\sigma$!).

As usual there exists a similar object $\rho'$ for the right invariant
functional $\psi$, i.e. $\rho'$ is an automorphism on $A$ such that
$\psi(a b) = \psi(b \rho'(a))$ for every $a,b \in A$.

Using the antipode, we can connect $\rho$ and $\rho'$ via the formula
$S \rho' = \rho S$. Furthermore, we have that $S^2$  
commutes with $\rho$ and $\rho'$.
The interplay between $\rho$,$\rho'$ and $\de$ is given by the following  
formulas:
$$\de \rho = (S^2 \od \rho)\de \hspace{1.5cm} \text{ and } \hspace{1.5cm}
\de \rho' = (\rho' \od S^{-2})\de.$$

\bigskip

It is also possible to introduce the modular function of our algebraic
quantum group. This is an invertible element $\sde$ in $M(A)$ such that
$(\vfi \od \io)(\de(a)(1 \ot b)) = \vfi(a) \, \sde  \, b $ for every $a,b
\in A$.

Concerning the right invariant functional, we have that
$(\io \od \psi)(\de(a)(b \ot 1)) = \psi(a) \, \sde^{-1} \, b$ for every
$a,b \in A$.

\medskip

This modular function is, like in the classical group case, a one
dimensional (generally unbounded)
corepresentation of our algebraic quantum group:
$$ \de(\sde) = \sde \ot \sde \hspace{1.5cm} \vep(\sde)= 1 \hspace{1.5cm}
S(\sde) = \sde^{-1} .$$

As in the classical case, we can relate the left invariant functional
to our right invariant functional via the modular function: we have
for every $a \in A$ that
$$ \vfi(S(a)) =  \vfi(a \sde) = \mu \, \vfi(\sde a) .$$
Not surprisingly, we have also that $\rho(\sde) = \rho'(\sde) = \mu^{-1}
\sde$.

\medskip

Another connection between $\rho$ and $\rho'$ is given by the equality
$\rho'(a) = \sde \rho(a) \sde^{-1}$ for all $a \in A$.

\bigskip

We have also a property which says, loosely speaking, that every element
of $A$ has compact support. This result was first proven by myself but
the (simpler) proof here is due to A. Van Daele.

\begin{proposition}
Consider $a_1,\ldots\!,a_n \in A$. Then there exists an element $c$ in
$A$ such that $c \, a_i = a_i \, c = a_i$ for every $i \in
\{1,\ldots\!,n\}$.
\end{proposition}
\begin{demo}
Define the following subspace $K$ of $A^{2n}$:
$$ K = \{ \,  (b a_1,\ldots\!,b a_n,a_1 b,\ldots\!,a_n b) \mid b \in A
\,\} \ .$$
We have to prove that $(a_1,\ldots\!,a_n,a_1,\ldots\!,a_n)$ belongs to $K$.

Therefore, choose a linear functional $\om$ on $A^{2n}$ which is 0 on
$K$. For every $i \in \{1,\ldots\!,2n\}$, we define $\om_i \in A'$ such
that $\om_i(x) = \om(0,\ldots\!, \stackrel{i}{x} , \ldots\! , 0)$ for
every $x \in A$.
It is clear that $\om(y) = \sum_{i=1}^{2n} \om_i(y_i)$ for every $y \in
A^{2n}$.

Take $d \in A$ such that $\vfi(d)=1$. We have for every $e \in A$ that
\begin{eqnarray*}
& & e \, [\,\sum_{i=1}^n (\om_i \od \io)(\de(d)(a_i \ot 1)) +
\sum_{i=1}^n (\om_{n+i} \od \io)((a_i \ot 1)\de(d)) \, ] \\ & &
\spat\!\!\! = \sum_{i=1}^n (\om_i \od \io)((1  \ot e) \de(d)(a_i \ot
1)) + \sum_{i=1}^n (\om_{n+i} \od \io)((a_i \ot 1)(1 \ot e)\de(d))
\end{eqnarray*}
which is 0 because  $\om=0$ on
$K$ and $(1 \ot e)\de(d)$ belongs to $A \od A$. So we get that
$$ \sum_{i=1}^n (\om_i \od \io)(\de(d)(a_i \ot 1)) + \sum_{i=1}^n
(\om_{n+i} \od \io)((a_i \ot 1)\de(d))  = 0  \ . $$
Applying $\vfi$ to this equality gives that
\begin{eqnarray*}
0 & = & \sum_{i=1}^n \om_i(\,(\io \od \vfi)(\de(d)(a_i \ot 1))\,) +
\sum_{i=1}^n \om_{n+i}(\,(\io \od \vfi)((a_i \ot 1)\de(d)) \,) \\ & = &
\sum_{i=1}^n \om_i(a_i) + \sum_{i=1}^n \om_{n+i}(a_i)
= \om(a_1,\ldots\!,a_n,a_1,\ldots\!,a_n) \ .
\end{eqnarray*}

From this all, we infer that $(a_1,\ldots\!,a_n,a_1,\ldots\!,a_n)$
belongs to $K$.
\end{demo}

\bigskip\medskip

In a last part of this section, we are going to say something about 
duality.

\begin{definition}
We define the subspace $\ah$ of $A'$ as follows:
$$\ah = \{\, \vfi a \mid a \in A \,\} = \{\, a \vfi \mid a \in A \,\} \ .$$
\end{definition}

The last equality results from the fact that $\vfi a = \rho(a) \vfi$ for
every $a \in A$.

\medskip

Because $\psi = \sde \vfi = \mu \, \vfi \sde$, we have also that 
$$\ah = \{\, \psi a \mid a \in A \,\} = \{\, a \psi \mid a \in A \,\}
\ .$$

\medskip

The faithfulness of $\vfi$ implies that $\ah$ separates $A$.

\bigskip

Like in the theory of Hopf-algebras, we want to turn $\ah$ into a
non-degenerate algebra:

\begin{itemize}
\item We know already that $(\io \od \om)\de(a)$ and $(\om \od
\io)\de(a)$ belong to $A$ for every $a \in A$ and $\om \in \ah$.
\item Choose $\om_1,\om_2 \in \ah$. Then there exist $a_1, a_2 \in A$
such that $\om_1 = \vfi a_1$ and $\om_2 = \vfi a_2$. Then it is not  
difficult to see that
$$\om_1( (\io \od \om_2)\de(x) ) = (\vfi \od \vfi)(\de(x)(a_1 \ot a_2))
= \om_2( (\om_1 \od \io)\de(x) ) $$
for every $x \in A$.
\end{itemize}

We have the following definition (see propositions 4.2 and 4.3 of \cite{VD1}).

\begin{definition} \label{def1.1}
We can turn $\ah$ into a non-degenerate algebra such that
$(\om_1 \om_2)(a) = \om_1( (\io \od \om_2)\de(a) ) = \om_2( (\om_1 \od
\io)\de(a))$ for every $\om_1,\om_2 \in \ah$ and $a \in A$
\end{definition}

\bigskip

In a next step,  a comultiplication is introduced on the level of $\ah$.
Because $\ah$ is a subspace of $A'$, we regard $\ah \od \ah$ in the
obvious way as a subspace of $(A \od A)'$. Therefore we can formulate the
next proposition (proposition 4.7 of \cite{VD1}).

\begin{proposition} \label{prop1.1}
There exists a unique non-degenerate homomorphism $\deh$ from $\ah$ into
$M(\ah \od \ah)$ such that we have for every $\om_1$, $\om_2$ in $\ah$ that
\begin{enumerate}
\item The element $\deh(\om_1)(1 \ot \om_2)$ belongs to $\ah \od \ah$ and
$[\deh(\om_1)(1 \ot \om_2)](x \ot y) = (\om_1 \od \om_2)((x \ot
1)\de(y))$ for every $x,y \in A$.
\item The element $(\om_1 \ot 1)\deh(\om_2)$ belongs to $\ah \od \ah$
and $[(\om_1 \ot 1)\deh(\om_2)](x \ot y) = (\om_1 \od \om_2)(\de(x)(1 \ot
y))$  for every $x,y \in A$.
\end{enumerate}
\end{proposition}

\bigskip

Next we define the linear functionals $\vfih$ and $\psih$ on
$\ah$ which will play the role of left and right invariant functionals on
$(\ah,\deh)$.

\begin{definition} \label{def1.2}
We define the linear functionals $\vfih$ and $\psih$ on $\ah$ such that
we have for every $a \in A$ that $\vfih(\psi \, a) = \vep(a)$ and
$\psih(a \, \vfi) = \vep(a)$.
\end{definition}

Now we can formulate a major result of \cite{VD1}.

\begin{theorem}
We have that $(\ah,\deh)$ is an algebraic quantum group with left Haar
functional $\vfih$ and right Haar functional $\psih$.
\end{theorem}

\begin{remark} \rm
The counit and antipode on $(\ah,\deh)$ are determined by the following
formulas:
\begin{itemize}
\item We have for every $a \in A$ that $\hat{\vep}(\vfi a) = \hat{\vep}(a
\vfi) = \vfi(a)$ and $\hat{\vep}(\psi a) = \hat{\vep}(a \psi) = \psi(a)$.
\item We have for every $\om \in \ah$ that $\hat{S}(\om) = \om \! \circ
\! S$.
\end{itemize}
The proof of these results can be found in \cite{VD1}.
\end{remark}

\section{Implementing the multiplier algebra of the dual as a space of
linear functionals}

Consider an algebraic quantum group $(A,\de)$ with left Haar
functional $\vfi$. We will use notations an conventions of the
previous section.

The non-degenerate algebra $\ah$ was introduced as a subset of $A'$.
In this section we will show that $M(\ah)$ can also be implemented as
a subset of $A'$.

\bigskip

\begin{definition} \label{def2}
We define left and right actions of $\ah$ on $A'$ as follows.
Consider $\om \in \ah$ and $\th \in A'$.
Then we define $\om \th$ and $\th \om$ in $A'$ such that
$(\om \th)(x) = \th( (\om \od \io)\de(x) )$ and
$(\th \om)(x) = \th( (\io \od \om)\de(x) )$ for every $x \in A$.
\end{definition}

If $\th$ belongs to $\ah$, these definitions of $\om \th$ and $\th
\om$ correspond to the ones given before.

\bigskip

The associativity of the product on $\ah$ will be essential in the
proof of the following lemma.

\begin{lemma} \label{lem4}
Consider $\om_1,\om_2 \in \ah$ and $\th \in A'$.
Then
\begin{enumerate}
\item $(\om_1 \om_2) \th = \om_1 (\om_2 \th)$
\item $\th (\om_1 \om_2) = (\th \om_1) \om_2$
\item $(\om_1 \th) \om_2 = \om_1 (\th \om_2)$
\end{enumerate}
\end{lemma}
\begin{demo}
Choose $x \in A$. We have for every $\om \in \ah$ that
\begin{eqnarray*}
\om( (\om_1 \om_2 \od \io) \de(x) ) & = & ((\om_1 \om_2) \om)(x)
= (\om_1 (\om_2 \om))(x) \\
& = & (\om_2 \om)( (\om_1 \od \io)\de(x) ) = \om(\, (\om_2 \od
\io)\de((\om_1 \od \io)\de(x)) \,)  \ .
\end{eqnarray*}
Because $\ah$ separates $A$, this implies that
$$ (\om_1 \om_2 \od \io) \de(x) = (\om_2 \od \io)\de((\om_1 \od
\io)\de(x)) \ .$$
Applying $\th$ to this equation gives in a similar way as above that
$((\om_1 \om_2) \th)(x) = (\om_1 (\om_2 \th))(x)$.

So we see that $(\om_1 \om_2) \th = \om_1(\om_2 \th)$.

\medskip

The other equalities are proven in an analogous way.
\end{demo}

\begin{lemma} \label{lem5}
Consider $\th \in A'$ such that $\th \om = 0$ for every $\om \in \ah$.
Then $\th = 0$.
\end{lemma}

This follows easily because we  can write every element of $A$ as a
sum of elements $(\io \od \vfi a)\de(b)$ with $a,b \in A$. Of course,
a similar separation result applies for left multiplication with
elements of $\ah$.

\bigskip

We will use the temporary notation
$$\tilde{A} = \{\, \th \in A' \mid \text{We have for every } \om \in \ah
\text{ that } \om \th \text{ and } \th \om \text{ belong to } \ah \,\} \ . $$
It is clear that $\ah$ is a subset of $\tilde{A}$.

\bigskip

Let $a$ be an element in $A$. For the purpose of clarity, we will use
the notations $\vfi [a] = \vfi a$ and $[a] \psi = a \psi$.

\medskip

\begin{lemma} \label{lem2}
Consider $a \in A$ and $\th \in A'$. Then the following equalities hold:
\begin{eqnarray*}
& & \bullet \ \ \th (\vfi a) = \vfi \, [S^{-1}\bigl((\io \od
\th)\de(S(a)) \bigr)] \hspace{2cm} \bullet \ \
\th(a \vfi) = \vfi \, [S\bigl((\io \od \th)\de(S^{-1}(a)) \bigr)]  \\
& & \bullet \ \ (\psi a) \th = \psi \, [S\bigl((\th \od
\io)\de(S^{-1}(a)) \bigr)] \hspace{1.965cm}  \bullet  \ \ (a \psi) \th
= [S^{-1}\bigl((\th \od \io)\de(S(a)) \bigr)]\, \psi
\end{eqnarray*}
\end{lemma}
\begin{demo}
We have for every $x \in A$ that
\begin{eqnarray*}
(\th (\vfi a))(x) & =  & \th( (\io \od \vfi a)\de(x) ) = \th(\,(\io
\od \vfi)((1 \ot a)\de(x)) \, )  \stackrel{(*)}{=} \th\bigl(\,S( \, (\io
\od \vfi)(\de(a) (1 \ot x))\,)\,\bigr) \\
& = & \vfi(\, (\th S \od \io)(\de(a)(1
\ot x)) \, ) = \vfi( \, (\th S \od \io)(\de(a)) \, x \, ) \ ,
\end{eqnarray*}
where we used equation \ref{eq1.1} in (*).
So $$ \th (\vfi a) = \vfi \, [ (\th S \od \io)\de(a)] = \vfi \,[
S^{-1}\bigl((\io \od \th)\de(S(a))\bigr)] \ , $$
where we used the equality $(S \od S)\de = \flip \de S$ in the last equation.

Using also equation \ref{eq1.2}, we find in a similar way the other results.
\end{demo}

\begin{lemma} \label{lem6}
Let $\th$ be a linear functional on $A$. Then $\th$ belongs to
$\tilde{A}$ $\Leftrightarrow$ We have for every $x \in A$ that
$(\th \od \io)\de(x)$ and $(\io \od \th)\de(x)$ belong to $A$.
\end{lemma}
\begin{demo}
\begin{description}
\item[\,\,$\Rightarrow$] \begin{minipage}[t]{15cm}
Choose $a \in A$. By assumption, there exists $b \in A$ such that $\th
(\vfi a) = \vfi b$.

This implies by the previous lemma
that $\vfi \, [S^{-1}\bigl((\io \od \th)\de(S(a))\bigr)] = \vfi b$. So
the faithfulness of $\vfi$ implies that
$S^{-1}\bigl((\io \od \th)\de(S(a)) \bigr) = b$,
so $(\io \od \th)\de(S(a)) = S(b)$ which belongs to $A$.

Using $\psi$ instead of $\vfi$, we prove in  a similar way that
$(\th \od \io)\de(S(a))$ belongs to $A$.
\end{minipage}
\item[\,\,$\Leftarrow$] \begin{minipage}[t]{15cm}
Choose $a \in A$. We have by assumption that $S^{-1}\bigl((\io \od
\th)\de(S(a))\bigr)$ belongs to $A$.

By the previous lemma, we get that $\th (\vfi a) = \vfi
\,[S^{-1}\bigl((\io \od \th)\de(S(a)) \bigr)]$
which belongs to $\ah$.

Similarly, we get
that $(\psi a) \th$ belongs to $\ah$.
\end{minipage}
\end{description}
\end{demo}

Because of lemma \ref{lem4}, we can define a mapping $F$ from
$\tilde{A}$ into $M(\ah)$ such that we have for every $\th \in
\tilde{A}$ and $\om \in A$ that $F(\th) \om = \th \om$ and $\om F(\th)
= \om \th$. Then we have that $F(\th) = \th$ for every $\th \in \ah$.

It is easy to check that $F$ is linear and injective (injectivity
follows from lemma \ref{lem5}). We will prove that $F$ is also
surjective.

\bigskip

Remember that we defined the non-zero linear functional $\vfih$ on
$\ah$ such that $\vfih(\psi a) = \vep(a)$ for every $a \in A$
(definition \ref{def1.2}).

\begin{lemma} \label{lem3}
Consider $a \in A$ and $\th \in A'$. Then $\vep( (\io \od \th)\de(a))=
\vep( (\th \od \io)\de(a) ) = \th(a)$.
\end{lemma}
\begin{demo}
There exist $c \in A$ such that $\vep(c)=1$.
Then
\begin{eqnarray*}
& & \vep( (\io \od \th)\de(a) ) =  \vep( (\io \od \th) \de(a) ) \,
\vep(c) = \vep( (\io \od \th)(\de(a)) \, c\,) \\
&  & \spat = \vep\bigl( (\io
\od \th)(\de(a)(c \ot 1) )\bigr)
= \th\bigl((\vep \od \io)(\de(a)(c \ot 1) ) \bigr)
=  \th( \vep(c) \, a) = \th(a) \ .
\end{eqnarray*}
The other equality is proven in a similar way.
\end{demo}

\begin{lemma} \label{lem1}
Consider $a \in A$ and $\om \in \tilde{A}$. Then $\vfih(\, (\psi a)\,
\om\,) = \om(S^{-1}(a))$.
\end{lemma}
\begin{demo}
By lemmas \ref{lem2} and \ref{lem6} , we know that $(\psi a) \, \om =
\psi \, [S\bigl((\om \od \io)\de(S^{-1}(a))\bigr)]$ and that
\newline $S\bigl((\om \od \io)\de(S^{-1}(a))\bigr)$ belongs to $A$. Therefore
$$
\vfih(\,(\psi a) \, \om\,)  =  \vep\bigl(\,S(\,(\om \od
\io)\de(S^{-1}(a))\,)\,\bigr) = \vep(\,(\om \od \io)\de(S^{-1}(a))\,)  $$
which by the previous lemma equals $\om(S^{-1}(a))$.
\end{demo}

This lemma suggests a solution for the proof of the surjectivity of
$F$. First, we  need another lemmas.

\begin{lemma} \label{lem7}
Consider $b \in A$ and $\om \in \ah$.
Then $\vfih(\,(\psi \sde S^2(b))\, \om\, ) = \vfih(\, \om (\psi b)\,)$.
\end{lemma}
\begin{demo} Take $a \in A$ such that $\om = \psi a$.
By the previous lemma we have that
\begin{eqnarray*}
\vfih(\,(\psi \sde S^2(b))(\psi a)\,) & = & (\psi a)(\,S^{-1}(\sde
S^2(b))\,) = (\psi a)(S(b) \sde^{-1}) \\
& = & \psi(a S(b) \sde^{-1}) = \vfi(a S(b)) = \psi(b S^{-1}(a)) \\
& = & (\psi b)(S^{-1}(a)) = \vfih(\,(\psi a)(\psi b)\,) \ ,
\end{eqnarray*}
where we used the previous lemma for a second time in the last equality.
\end{demo}

\begin{lemma}
Consider $b,x \in A$. Then
$\psi \, [S\bigl((\io \od \psi b)\de(x)\bigr)]  = (\psi \sde S^2(b))(\psi
S(x))$.
\end{lemma}
\begin{demo}
By lemma \ref{lem2}, we know that
$$(\psi \sde S^2(b))(\psi S(x)) = \psi \, [S\bigl(\,(\psi S(x) \od
\io)\de(S^{-1}(\sde S^2(b)))\,\bigr)] =  \psi \, [S\bigl((\psi S(x) \od
\io)\de(S(b) \sde^{-1})\bigr)] \ .$$
But we have that
\begin{eqnarray*}
& & S\bigl((\psi S(x) \od \io)\de(S(b) \sde^{-1})\bigr) = S\bigl((\psi
\od \io)((S(x) \ot 1)\de(S(b) \sde^{-1})) \bigr) \\ & & \spat
\stackrel{(*)}{=} (\psi \od \io)(\de(S(x)) (S(b) \sde^{-1} \ot 1)) =
(\sde^{-1} \psi \od \io)(\de(S(x))(S(b) \ot 1)) \\ & & \spat =  (\vfi
\od \io)(\,\flip(S \od S)((1 \ot b)\de(x))\,)
=  S\bigl((\io \od \vfi S)((1 \ot b)\de(x))\bigr) \\
& & \spat = S\bigl((\io \od \psi)(( 1 \ot b)\de(x))\bigr)
= S\bigl((\io \od \psi b)\de(x)\bigr)
\end{eqnarray*}
where we used equation \ref{eq1.2} in equality (*). Combining these two
results, the lemma follows.
\end{demo}

\medskip

We get now to a proposition of which the proof is rather
straightforward.

\begin{proposition}
The mapping $F$ is a isomorphism of vector spaces from $\tilde{A}$ to
$M(\ah)$.
\end{proposition}
\begin{demo}
Choose $T \in M(\ah)$. Define the element $\th \in A'$ such that
$\th(x) = \vfih(\,(\psi S(x))\,T\,)$ for every $x \in A$ (This is
suggested by lemma \ref{lem1}). Choose $b \in A$.
\begin{enumerate}
\item First we prove that $(\psi b) \th = (\psi b) \, T$

Choose $x \in A$. Using definition \ref{def2}, we get that
$$ ((\psi b) \th)(x) = \th( (\psi b \od  \io)\de(x) ) =  \vfih(\,(\psi
\,[S\bigl((\psi b \od \io)\de(x)\bigr)]) \, T \,) \ .$$
By lemma \ref{lem2}, we have that $\psi \, [S\bigl((\psi  b \od
\io)\de(x)\bigr)] = (\psi S(x)) (\psi b)$.
Therefore, the associativity of the product in $M(\ah)$ and lemma  
\ref{lem1} imply that
$$ ((\psi b) \th)(x) = \vfih(\,[(\psi S(x)) (\psi b)] \, T\,)
= \vfih(\, (\psi S(x)) \, [ (\psi b) T ]\,)
= ((\psi b) \, T)(x) \ . $$

So we see that $(\psi b) \th = (\psi b) \, T$.

\item Next we prove that $\th (\psi b) = T \, (\psi b)$.

Choose $x \in A$. Then definition \ref{def2} implies that
$$ (\th (\psi b))(x) = \th( (\io \od \psi b)\de(x) )
= \vfih(\,(\psi \, [S\bigl((\io \od \psi b)\de(x)\bigr)]) \,\, T\,) \ .$$
The previous lemma implies that $\psi\,[S\bigl((\io \od \psi
b)\de(x)\bigr)] = (\psi \sde S^2(b))(\psi S(x))$.

So we see that
$$(\th (\psi b))(x) = \vfih( \,[(\psi \sde S^2(b))(\psi S(x))] \,T\,)
= \vfih( \, (\psi \sde S^2(b)) \,[ (\psi S(x)) T ]\,) \ .$$
Using lemma \ref{lem7} and \ref{lem1}, this equation implies that
$$(\th (\psi b))(x) = \vfih(\, [(\psi S(x)) T] \,(\psi b) \,)
= \vfih(\,(\psi S(x))\, [T (\psi b)] \, ) = (T \, (\psi b))(x) \ .$$

So we arrive at the conclusion that $\th (\psi b) = T \, (\psi b)$

\end{enumerate}

Looking at the definition of $\tilde{A}$ and $F$, we get from these two
equations easily that $\th$ belongs to $\tilde{A}$
and that $F(\th) = T$.
\end{demo}

\medskip

Now we will use the mapping $F$ to turn $\tilde{A}$ into an algebra.
We define on $\tilde{A}$ a product operation such that
$F(\th_1 \th_2) = F(\th_1) F(\th_2)$ for every $\th_1,\th_2 \in \tilde{A}$.

\medskip

\begin{proposition}
We have for every $\th_1,\th_2 \in \tilde{A}$ and $x \in A$ that
$(\th_1 \th_2)(x) = \th_1( (\io \od \th_2)\de(x) ) = \th_2( (\th_1 \od
\io)\de(x) )$.
\end{proposition}
\begin{demo}
 Choose $a \in A$. Using lemma \ref{lem2} twice, we get that
\begin{eqnarray*}
F(\th_1 \th_2) \, (\vfi a) & = & (F(\th_1) F(\th_2))\, (\vfi a) =
F(\th_1) \, (F(\th_2) (\vfi a))  \\
& = & F(\th_1)\, ( \th_2 (\vfi a))
=  F(\th_1) \, \, \vfi \, [S^{-1}\bigl((\io \od \th_2)\de(S(a))\bigr)]   \\
& = &\th_1 \, \, \vfi \, [S^{-1}\bigl((\io \od \th_2)\de(S(a))\bigr)]
\\ & = & \vfi [S^{-1}\bigl(\,(\io \od \th_1)\de( \,(\io \od
\th_2)\de(S(a))\,)\,\bigr)] \ .
\end{eqnarray*}

On the other hand, lemma \ref{lem2} implies also that
$$F(\th_1 \th_2) \, (\vfi a) =
(\th_1 \th_2) \, (\vfi a) = \vfi \, [S^{-1}\bigl((\io \od \th_1
\th_2)\de(S(a))\bigr)] \ . $$
Comparing these two different expressions for $F(\th_1 \th_2)$ and
remembering that $\vfi$ is faithful, we get that
$$S^{-1} \bigl(\,(\io \od \th_1)\de(\,(\io \od \th_2)\de(S(a))\,)\,
\bigr) = S^{-1}\bigl((\io \od \th_1 \th_2)\de(S(a))\bigr) \ ,$$
which implies that
$$(\io \od \th_1)\de((\io \od \th_2)\de(S(a)))
= (\io \od \th_1 \th_2)\de(S(a)) \ .$$
Applying $\vep$ to this equation and using lemma \ref{lem3}, we see that
$\th_1((\io \od \th_2)\de(S(a)))
= (\th_1 \th_2)(S(a))$.

\medskip

The other equality is proven in a similar way.
\end{demo}

\bigskip

Until now, we used $M(\ah)$ as an abstract object. From now on, we will
use $\tilde{A}$ as a concrete realization of $M(\ah)$ and want to forget
about the mapping $F$. This will be formulated in the next theorem.

\begin{theorem} \label{thm1}
As a vector space, $M(\ah)$ is equal to
$$\{\, \th \in A' \mid \text{ We have for every } a \in A
\text{ that } (\th \od \io)\de(a) \text{ and } (\io \od \th)\de(a) 
\text{ belong to } A \,\} \ . $$
The product operation on $M(\ah)$ is defined in such a way
that  $(\th_1 \th_2)(a) = \th_1( (\io \od \th_2)\de(a) )$ \newline $=
\th_2( (\th_1 \od \io)\de(a) )$
for every $\th_1,\th_2 \in M(\ah)$ and $a \in A$.
Furthermore, $\vep$ is the unit of $M(\ah)$.
\end{theorem}

The last statement of the proposition is easy to check. The others have
been proven in the previous part of the section.

\medskip

\begin{remark} \rm
Consider $a \in A$ and $\om \in A$. Using the previous theorem, it is
easy to check that $\om a$ and $a \om$ belong to $M(\ah)$.
\end{remark}

\bigskip

From this section, we only need to remember the previous theorem and
the following results. The rest of this section was intended to prove
these results.

\medskip

We can also extend the left and right actions of $\ah$ on $A'$
to left and right actions of $M(\ah)$ on $A'$ in the obvious way.

\begin{definition} \label{def1}
Consider $\om \in M(\ah)$ and $\th \in A'$. We define $\om \th$ and $\th
\om$ in $A'$ such that $(\om \th)(a) = \th((\om \od \io)\de(a))$ and
$(\th \om)(a) = \th((\io \od \om)\de(a))$ for every $a \in A$.
\end{definition}

This definition implies for instance that $((\om b) \th)(x) = (\om \od
\th)( (b \ot 1) \de(a) )$ for every $\om,\th \in A'$ and $a,b \in A$.
Another consequence is the equality $\th \vep = \vep \th = \th$ for every
$\th \in A'$.

\medskip

The proof of the following proposition is the same as the proof of  lemma
\ref{lem4}.

\begin{proposition} \label{prop1}
Consider $\om_1,\om_2 \in M(\ah)$ and $\th \in A'$.
Then
\begin{enumerate}
\item $(\om_1 \om_2) \th = \om_1 (\om_2 \th)$
\item $\th (\om_1 \om_2) = (\th \om_1) \om_2$
\item $(\om_1 \th) \om_2 = \om_1 (\th \om_2)$
\end{enumerate}
\end{proposition}

\medskip

Remember that we have also the following result.

\begin{proposition} \label{prop2}
Consider $\th \in A'$. Then $\th$ belongs to $M(\ah)$
$\Leftrightarrow$ We have for every $\om \in \ah$ that
$\th \om$ and $\om \th$ belong to $\ah$.
\end{proposition}

\section{A nice implementation of the comultiplication on the dual}

Consider an algebraic quantum group $(A,\de)$ with a left Haar
functional $\vfi$ (we will use the notations for counit, antipode, ...
as in the previous sections). The results of the previous section will
be used to get  nice formulas for the comultiplication, counit and
antipode on $\ah$.

\bigskip

As can be expected, we put on $A \od A$ a comultiplication $\de$ such
that $\de(a \ot b) = (\io \od \flip \od \io)(\de(a) \ot \de(b))$ for
every $a,b \in A$ (The use of the same symbol for the comultiplication
on A as for the comultiplication on $A \od A$ should not cause any
confusion because we always know on which elements they work).

It is not so hard to check that $(A \od A,\de)$ is again an algebraic
quantum group with counit $\vep \od \vep$, antipode $S \od S$ and left
Haar functional $\vfi \od \vfi$. Moreover, $(A \od A)\hoed$ will be
equal to $\ah \od \ah$.

By the results of the previous section, we know that $M(\ah \od \ah)$
is a subset of $(A \od A)'$ and this will be essential to the next
proposition.

\begin{lemma}
Let $a,b$ be elements of $A$ and $\om,\th$ elements in $A'$. Then
$(\io \od \io \od  \om \ot \th)\de(a \ot b) =$ \newline $(\io \od
\om)\de(a)\ot (\io \od \th)\de(b)$ and
$(\om \ot \th \od \io \od \io)\de(a \ot b) = (\om \od \io)\de(a)
\ot (\th \od \io)\de(b)$.
\end{lemma}
\begin{demo}
Choose $c,d \in A$. Then we have by definition that
$$(c \ot d) \, (\io \od \io \od \om \ot \th)\de(a \ot b)
= (\io \od \io \od \om \ot \th)((c \ot d \ot 1 \ot 1)\de(a \ot b)) \ . $$
We have also that
$$(c \ot d \ot 1 \ot 1)\de(a \ot b)
= (c \ot d \ot 1 \ot 1) (\io \od \flip \od \io)(\de(a) \ot \de(b))
= (\io \od \flip \od \io)( (c \ot 1)\de(a) \ot (d \ot 1)\de(b) )  \ . $$
Notice that $(c \ot 1)\de(a)$ and $(d \ot 1)\de(b)$ belong to $A \od A$.
So we get that
\begin{eqnarray*}
& & (c \ot d) \, (\io \od \io \od \om \ot \th)\de(a \ot b)
= (\io \od \io \od \om \ot \th)\bigl((\io \od \flip \od \io)( (c \ot
1)\de(a) \ot (d \ot 1)\de(b) )\bigr) \\
& & \spat = (\io \od  \om)((c \ot 1)\de(a)) \ot (\io \od \th)((d \ot
1)\de(b))
= (c \ot d) \, ((\io \od \om)\de(a) \ot (\io \od \th)\de(b)) \ .
\end{eqnarray*}
So we get that
$(\io \od \io \od \om \ot \th)\de(a \ot b) = (\io \od \om)\de(a)
\ot (\io \od \th)\de(b)$. The other equality is proven in a similar way.
\end{demo}

\begin{proposition}
Consider $\om \in \ah$. Then we have for every $a,b \in A$ that
$\hat{\de}(\om)(a \ot b) = \om(a b)$.
\end{proposition}
\begin{demo}
Define $\th \in (A \od A)'$ such that $\th(a \ot b) = \om(a b)$ for every
$a,b \in A$.

Take $\eta_1,\eta_2 \in \ah$. Choose $x,y \in A$. We know that $(\io
\od \eta_2)\de(y)$ belongs to $A$. By definition \ref{def1}, we get
that
\begin{eqnarray*}
& & (\th \, (\vep \ot \eta_2))(x \ot y)
= \th((\io \od \io \od \vep \ot \eta_2)\de(x \ot y)) \\
& & \spat = \th( (\io \od \vep)\de(x) \ot (\io \od \eta_2)\de(y)\,)
= \th( x \ot (\io \od \eta_2)\de(y)) \\
& & \spat = \om(x \, (\io \od \eta_2)\de(y))
= \om((\io \od \eta_2)((x \ot 1)\de(y))
= (\om \od \eta_2)((x \ot 1)\de(y)) \ .
\end{eqnarray*}
Hence, proposition \ref{prop1.1} implies that
$$(\th \, (\vep \ot \eta_2))(x \ot y) = (\deh(\om)(\vep \ot \eta_2))(x
\ot y) \ . $$
Therefore, we get that $\th \, (\vep \ot \eta_2) = \deh(\om)(\vep \ot
\eta_2)$. Multiplying this equation with $\eta_1 \ot \vep$ to the right
and using proposition \ref{prop1}, we see that
$$\th \, (\eta_1 \ot \eta_2) = \deh(\om)(\eta_1 \ot \eta_2) \ .$$
By lemma \ref{lem5}, this implies that $\th = \deh(\om)$.
\end{demo}

We can even do better:

\begin{proposition}
Consider $\om \in M(\ah)$. Then we have for every $a,b \in A$ that
$\hat{\de}(\om)(a \ot b) = \om(a b)$.
\end{proposition}
\begin{demo}
Define $\th \in (A \od A)'$  such that $\th(a \ot b) = \om(a b)$ for every
$a,b \in A$.

Choose $y \in A$ and $\om_1,\om_2 \in \ah$.
Then we have for $a,b \in A$ that
$$
(\th \, \deh(\vfi y))(a \ot b) =
\th(\, (\io \od \io \od \deh(\vfi y))\de(a \ot b) \,) \ .
\text{\ \ \ \ \ \ \ \ (a)}
$$
Notice that $(1 \ot 1 \ot y)\de(a)_{13} \de(b)_{23}$ belongs to $A \od A
\od A$.
We have for every $c,d \in A$ that
\begin{eqnarray*}
[(\io \od \io \od \deh(\vfi y))\de(a \ot b)] \, (c \ot d) & = & (\io
\od \io \od \deh(\vfi y))(\de(a \ot b) (c \ot d \ot 1 \ot 1)) \\ & = &
(\io \od \io \od \deh(\vfi y))(\,(\io \od \flip \od \io)(\de(a)(c \ot
1) \ot \de(b)(d \ot 1)) \,) \\ & \stackrel{(*)}{=} & (\io \od \io \od
\vfi) (\, (1 \ot 1 \ot y)(\de(a)(c \ot 1))_{13} (\de(b)(d \ot 1))_{23}
\,) \\ & = & (\io \od \io \od \vfi)((1 \ot 1 \ot y)\de(a)_{13}
\de(b)_{23})\, (c \ot d) \ .
\end{eqnarray*}
Here we used the previous lemma in equality (*).

So we see that
$$(\io \od \io \od \deh(\vfi y))\de(a \ot b) =
(\io \od \io \od \vfi)((1 \ot 1 \ot y)\de(a)_{13} \de(b)_{23}) \ . $$
Plugging this equality in equality (a) gives
$$
(\th \, \deh(\vfi y))(a \ot b)
= \vfi(\, (\th \od \io)((1 \ot 1 \ot y)\de(a)_{13} \de(b)_{23})\,) \ .
\text{\ \ \ \ \ \ (b)}$$
The definition of $\th$ implies for every $e \in A$ that
\begin{eqnarray*}
(\th \od \io)((1 \ot 1 \ot y)\de(a)_{13} \de(b)_{23}) \,\, e
& = & (\th \od \io)( [(1 \ot y)\de(a)]_{13} [\de(b)(1 \ot e)]_{23}) \\
& = & (\om \od \io \od \io)( [(1 \ot y)\de(a)]\,[\de(b)(1 \ot e)] ) \\
& = & (\om \ot \io \od \io)( (1 \ot y)\de(a b)) \,\, e \ .
\end{eqnarray*}
So we get that $(\th \od \io)((1 \ot 1 \ot y)\de(a)_{13} \de(b)_{23})
= (\om \od \io)( (1 \ot y)\de(a b))$.
Using this in equality (b), gives us that
$$ (\th \, \deh(\vfi y))(a \ot b)
= (\om \od \vfi)( (1 \ot y)\de(a b) )
= (\om \, (\vfi y))(ab) \ .$$
Therefore, the previous lemma implies that
$$ (\th \, \deh(\vfi y))(a \ot b) = \deh(\om \, (\vfi y))(a \ot b) \ .$$

From this all we get that
$$\th \, \deh(\vfi y) = \deh(\om \,(\vfi y)) = \deh(\om) \deh(\vfi y) \ .$$
Multiplying this equation from the right by $\om_1 \ot \om_2$ and
using proposition \ref{prop1}, results in the equality
$$\th \, [\deh(\vfi y)(\om_1 \ot \om_2)] = \deh(\om)
\, [\deh(\vfi y)(\om_1 \ot \om_2)] \ . $$

Using the non-degeneracy of $\deh$ and lemma \ref{lem5}, we infer from
this all that $\deh(\om) = \th$.
\end{demo}

\medskip

This proposition implies easily the following result.

\begin{result}
Consider $\om \in M(\ah)$ and $\om_1,\om_2 \in M(\ah)$. Then
$[\deh(\om)(\om_1 \ot \om_2)](x \ot y) =$ \newline $\om\bigl(\,(\io
\od \om_1)\de(x) \,\, (\io \od \om_2)\de(y)\, \bigr)$ for every $x,y \in A$.
\end{result}
\begin{demo}
Remember that $(\io \od \om_1)\de(x)$ and $(\io \od \om_2)\de(y)$
belong to $A$. So we get that
\begin{eqnarray*}
& & [\deh(\om)(\om_1 \ot \om_2)](x \ot y)
= \deh(\om)( (\io \od \io \od \om_1 \ot \om_2)\de(x \ot y) ) \\
& & \spat = \deh(\om)( (\io \od \om_1)\de(x) \ot (\io \od \om_2)\de(y))
= \om((\io \od \om_1)\de(x) \,\, (\io \od \om_2)\de(y)) \ ,
\end{eqnarray*}
where we used the previous proposition in the last equality.
\end{demo}

\begin{remark} \rm \label{rem2.1}
A special case of the foregoing result is the case where $\om_1=\vep$ or
$\om_2=\vep$. This gives  rise to the following equalities:

\medskip

Consider $\om_1,\om_2 \in M(\ah)$. Then we have for every $x,y \in A$
that $[\deh(\om_1)(\vep \ot \om_2)](x \ot y) =$ \newline $(\om_1 \od
\om_2)((x \ot 1)\de(y))$ and $[\deh(\om_1)(\om_2 \ot \vep )](x \ot y) =
(\om_1 \od \om_2)(\de(x)(y  \ot 1))$.
\end{remark}

\medskip

Of course, a similar calculation as in the proof  gives rise to
formulas for $(\om_1 \ot \om_2)\deh(\om)$ for every $\om,\om_1,\om_2
\in M(\ah)$.

\bigskip\medskip

In the last propositions of this section, we prove the usual formulas
for the counit and antipode on $M(\ah)$.

\begin{proposition}
We have for every $\om \in A'$ and $a \in A$ that
$\hat{\vep}(\om a) = \om(a)$ and $\hat{\vep}(a \om) = \om(a)$.
\end{proposition}
\begin{demo}
Take $b \in A$ such that $\vfi(b)=1$. Then $\hat{\vep}(\vfi b) = \vfi(b)
= 1$.

By lemma \ref{lem2}, we have that $(\om a) (\vfi b) = \vfi \,
[S^{-1}\bigl((\io \od \om a)\de(S(b))\bigr)]$. This implies that
\begin{eqnarray*}
\hat{\vep}(\om a) & = & \hat{\vep}((\om a)(\vfi b))
= \hat{\vep}(\,\vfi \, [S^{-1}\bigl((\io \od \om a)\de(S(b))\bigr)]\,)
= \vfi\bigl(\,S^{-1}(\,(\io \od \om a)\de(S(b))\,)\,\bigr) \\
& = & (\vfi S^{-1})\bigl(\,(\io \od \om)(\,(1 \ot a)\de(S(b))\,)\,\bigr)
= \om\bigl(\,(\vfi S^{-1} \od \io)(\,(1 \ot a)\de(S(b))\,)\, \bigr) \ .
\end{eqnarray*}
Using the right invariance of $\vfi S^{-1}$, this implies that
$$\hat{\vep}(\om a) = \om( \, a \,\, (\vfi S^{-1})(S(b))\,) = \om(a \,
\vfi(b)) = \om(a) \ . $$
The other equality is proven in a similar way.
\end{demo}

\bigskip

\begin{proposition}
Consider $\om \in M(\ah)$. Then $\hat{S}(\om)  = \om \! \circ \! S$
and $\hat{S}^{-1}(\om) = \om \! \circ \! S^{-1}$.
\end{proposition}
\begin{demo}
It is not difficult to check that $\hat{S}(\vfi c) = (\vfi c) \! \circ
\! S = S^{-1}(c) \psi$ for every $c \in A$.

Take $b \in A$.

Using lemma \ref{lem2} once again, we see that
$\om \, (\vfi b) = \vfi \, [S^{-1}\bigl((\io \od \om)\de(S(b))\bigr)]$.
By the remark in the beginning of this proof, we see that
$$\hat{S}(\vfi b) \, \hat{S}(\om) = \hat{S}(\om \, (\vfi b))=
[S^{-2}\bigl((\io \od \om)\de(S(b))\bigr)] \, \psi \ .$$

On the other hand, $\hat{S}(\vfi b) \, \hat{S}(\om) = (S^{-1}(b) \psi)
\, \hat{S}(\om)$ which by lemma \ref{lem2} implies that
$$\hat{S}(\vfi b) \, \hat{S}(\om) = [S^{-1}\bigl((\hat{S}(\om) \od
\io)\de(b)\bigr)] \, \psi \ .$$
Comparing these two expression for $\hat{S}(\vfi b) \, \hat{S}(\om)$ and
using the faithfulness of $\psi$, we arrive at the conclusion that
$$S^{-2}\bigl((\io \od \om)\de(S(b))\bigr) = S^{-1}\bigl((\hat{S}(\om)
\od \io )\de(b)\bigr) \ . $$
If we apply $\vep$ to this equation and use lemma \ref{lem3}, we see that
$\om(S(b)) = \hat{S}(\om)(b)$.

So we see that $\hat{S}(\om) = \om \! \circ \! S$.

The result concerning $\hat{S}^{-1}$ follows immediately from the
result concerning $\hat{S}$.
\end{demo}

\section{Slicing with certain functionals}

In this section, we will consider an algebraic quantum group $(A,\de)$
with left Haar functional $\vfi$.

\medskip

Let $B$ be a non-degenerate algebra and  $V$ an element of $M(A \od B)$.
In this section, we want to define a left multiplier $(\om \od \io)(V)$
for certain functionals $\om$ on $A$ and prove some familiar calculation
rules about this kind of slicings.

\medskip

\begin{notation}
We define the set $A^\circ$ as the vector space $\langle \, a \om
\mid a \in A \, \rangle$. So $A^\circ$ is a subspace of $A'$.
\end{notation}

By  theorem \ref{thm1}, it is clear that $A^\circ$ is a subspace of
$M(\ah)$. It is not so difficult to check that $A^\circ$ is a
subalgebra of $M(\ah)$. We have also immediately that $\ah$ is a
subset of $A^\circ$.

\medskip

We will use these functionals to make slicings. First we need an easy but
useful result.

\begin{lemma}
Consider a non-degenerate algebra $B$ and $V \in M(A \od B)$. Let
$a_1,\ldots\!,a_n \in A$, and $\om_1,\ldots\!,\om_n \in A'$ such that
$\sum_{i=1}^n a_i \, \om_i = 0$. Then $\sum_{i=1}^n (\om_i \od \io)(
V(a_i \ot x)) = 0$ for every $x \in B$.
\end{lemma}
\begin{demo}
Choose $x \in B$. There exists $e \in A$ such that $e a_i  = a_i$ for
every $i \in \{1,\ldots\!,n\}$. Then we have that $V (e \ot x)$ belongs
to $A \od B$. This implies that
\begin{eqnarray*}
\sum_{i=1}^n  (\om_i \od \io)(V(a_i \ot x)) & = & \sum_{i=1}^n (\om_i \od
\io)(V(e a_i \ot x)) = \sum_{i=1}^n (a_i \om_i \od \io)(V(e  \ot x)) \\
& = & ((\sum_{i=1}^n a_i \, \om_i) \od \io)(V(e \ot x)) = 0 \ .
\end{eqnarray*}
\end{demo}

\medskip

\begin{proposition}
Consider a non-degenerate algebra $B$ and $V \in M(A \od B)$. There
exists a unique linear map $G$ from $A^\circ$ into the set of left
multipliers on $B$ such that $G(a \om)\,x  = (\om \od \io)(V(a \ot
x))$  for every $\om \in A$, $a \in A$ and $x \in B$.

For every $\th \in A^\circ$, we put $(\th \od \io)(V) = G(\th)$, so
$(\th \od \io)(V)$ is a left multiplier on $B$.
\end{proposition}
\begin{demo}
Using the definition of left multipliers on $B$ it is not difficult to
check for every $a \in A$ and $\om \in A$ the existence of a unique
left multiplier $T(a,\om)$ on $B$ such that $T(a,\om)\,x = (\om \od
\io)(V(a \ot x))$  for every $x \in B$.

Let $a_1,\ldots\!,a_m , b_1,\ldots\!,b_n \in A$, $\om_1,\ldots\!,\om_m,
\th_1,\ldots\!, \th_n \in A'$ such that $\sum_{i=1}^m a_i \, \om_i =
\sum_{j=1}^n b_j \, \th_j $. Then the previous lemma implies that
$\sum_{i=1}^m T(a_i,\om_i) = \sum_{j=1}^n T(c_j,\th_j)$.

So we can define a mapping $G$ from $A^\circ$ into the set of left
multipliers on $B$ such that $G(\sum_{i=1}^m a_i \, \om_i ) =
\sum_{i=1}^m T(a_i,\om_i)$ for every $a_1,\ldots\!,a_m  \in A$ and
$\om_1,\ldots\!,\om_n \in A'$. It is easy to check that $G$ satisfies
the requirements of the proposition.
\end{demo}

It is easy to see that this notation $(\th \od \io)(V)$ is consistent
with the usual notation in the case that $V$ is an element of $A \od
B$.

\bigskip

\begin{remark}
Consider a non-degenerate algebra $B$. We only need to remember the
following obvious rules:
\begin{itemize}
\item The mapping $A^\circ \times M(A \od B) \rightarrow \text{ the set
of left multipliers of $B$} : (V,\om) \mapsto (\om \od \io)(V)$ is bilinear.
\item Consider $V \in M(A \od B)$, $a \in A$ and $\om \in A'$. Then
$(a \om \od \io)(V) \, x = (\om \od \io)(V(a \ot x))$ for every $x \in B$.
\end{itemize}
\end{remark}

\bigskip

The usual separation results stay valid.

\begin{result} \label{res3.3}
Consider a non-degenerate algebra $B$ and $V \in M(A \ot B)$. Then
$V = 0$ $\Leftrightarrow$ $(\om \od \io)(V) = 0$ for every $\om \in \ah$.
\end{result}
\begin{demo}
Suppose that $(\om \od \io)(V) = 0$ for every $\om \in \ah$.

Take $a \in A$, $b \in B$. Choose $\eta \in B'$.
Then we have for every $c \in A$ that
\begin{eqnarray*}
& & \vfi((\io \od \eta)(V (a \ot b))\,c)
= \vfi\bigl((\io \od \eta)(V(a c \ot b))\bigr) \\
& & \spat = \eta\bigl((\vfi \od \io)(V(a c \ot b))\bigr)
= \eta((a c \vfi \od \io)(V) \, b ) = 0 \,
\end{eqnarray*}
The faithfulness of $\vfi$ implies that $(\io \od \eta)(V(a \ot b)) = 0$.
Hence, $V (a \ot b)= 0$.

So we see that $V = 0$.
\end{demo}

We will be mainly interested in the cases where $(\om \ot \io)(V)$
becomes an element in $M(B)$. The following result deals with a natural
case. The proof of this result is straightforward and will therefore be
left out.

\begin{result} \label{res3.1}
Consider a non-degenerate algebra $B$ and $V \in M(A \od B)$. Let $a,b$
be elements in $A$ and $\om$ an element in $A'$. Then $(a \om b \od
\io)(V)$ is an element of $M(B)$ and we have for every $x \in B$ that $x
\, (a \om b \od \io)(V) = (\om \od \io)( (b \ot x) V (a \ot 1))$ and $(a
\om b \od \io)(V) \, x = (\om \od \io)( (b \ot 1) V (a \ot x))$.
\end{result}

\begin{corollary}
Consider a non-degenerate algebra $B$ and $V \in M(A \od B)$. Then we
have for every $\om \in \ah$ that $(\om \od \io)(V)$ belongs to $M(B)$.
\end{corollary}

\medskip

\begin{result}
Consider a non-degenerate algebra $B$ and $V \in M(A \od B)$. Then we
have for every $a \in A$ and $b \in B$ that
$(a \vfi \od \io)(V) \, b = (\vfi \od \io)(V(a \ot b))$ and
$b \, (\vfi a \od \io)(V) = (\vfi \od \io)((a \ot b) V)$.
\end{result}
\begin{demo}
The first result is true by definition. Let us turn to the second one.

There exist $e \in A$ such that $a =  e a$, then $\vfi a = \vfi e a =
\rho(e) \vfi a$. By result \ref{res3.1} , this implies for every $b \in
B$ that
\begin{eqnarray*}
b \, (\vfi a \od \io)(V) & = & b \, (\rho(e) \vfi a \od \io)(V)
= (\vfi \od \io)( (a \ot b) V (\rho(e) \ot 1)) \\
& = & (\vfi \od \io)( (e a \ot b) V ) = (\vfi \od \io)((a \ot b)V) \ .
\end{eqnarray*}
\end{demo}

\medskip

In the next result, we prove a natural calculation rule.

\begin{result} \label{res3.2}
Consider non-degenerate algebras $B$,$C$ and a non-degenerate
homomorphism $\pi$ from $B$ into $M(C)$. Let $V$ be an element in $M(A
\ot B)$ and $\om \in A^\circ$. Then
$\pi((\om \od \io)(V)) = (\om \od \io)((\io \od \pi)(V))$.
\end{result}
\begin{demo}
Choose $a \in A$ and $\th \in A'$.
Then we have for every $x \in B$ that
\begin{eqnarray*}
& & \pi((a \th  \od \io)(V)) \, \pi(x) y
= \pi(\,(a \th  \od \io)(V)\, x \,) \,y
= \pi(\,(\th \od \io)(V(a \ot x))\,) \, y\\
& & \spat = (\th \od \io)(\,(\io \od \pi)(V(a \ot x))\, (1 \ot y) \,) =
(\th \od \io)(\, (\io \od \pi)(V)(a \ot \pi(x)y)\,)\\
& & \spat = (a \th \od \io)((\io \od \pi)(V)) \, \pi(x) y
\end{eqnarray*}
Therefore, the non-degeneracy of $\pi$ implies that
$\pi((a \th \od \io)(V)) = (a \th \od \io)((\io \od \pi)(V))$.
\end{demo}

\section{Corepresentations of algebraic quantum groups}

Consider an algebraic quantum group $(A,\de)$ with a left Haar
functional $\vfi$. We will introduce the notion of a corepresentation
of $(A,\de)$ and prove some familiar results about them.

\medskip

We start of with the usual definition of a corepresentation.

\begin{definition}
Consider a non-degenerate algebra $B$. A corepresentation of $(A,\de)$ on
$B$ is by definition an element $V$ in $M(A \od B)$ such that
$(\de \od \io)(V) = V_{13} V_{23}$.

We call $V$ non-degenerate $\Leftrightarrow$ $V$ is invertible in $M(A
\od B)$.
\end{definition}

\medskip

Remember that $A \od A$ is again an algebraic quantum group  with $(A
\od A)\hoed = \ah \od \ah$. It is also easy to check that $A^\circ \od
A^\circ \subseteq (A \od A)^\circ$.

\medskip

Let $B$ be a non-degenerate algebra and $V \in M(A \od A \od B)$. The
remark above implies that we can define the left multiplier $(\om \od
\io)(V)$ on $B$ for every $\om \in A^\circ \od A^\circ$ (using the
definitions of the previous section).

\begin{lemma}  \label{lem4.1}
Consider a non-degenerate algebra $B$ and $V \in M(A \od B)$. Then we
have for every \newline $\om_1,\om_2 \in A^\circ$ that $(\om_1 \om_2
\od \io)(V) = (\om_1 \od \om_2 \od \io)((\de \ot \io)(V))$.
\end{lemma}
\begin{demo}
Take $\th_1,\th_2 \in A'$ and $a_1,a_2 \in A$. Then there exist
$p_1,\ldots\!,p_n , q_1,\ldots,q_n , r_1,\ldots\!,r_n \in A$ such that
$a_1 \ot a_2 = \sum_{i=1}^n \de(r_i)(p_i \ot q_i)$.

Define for every $i \in \{1,\ldots\!,n\}$ the element $\eta_i \in A'$
such that $\eta_i(x) = (\th_1 \od \th_2)(\de(x)(p_i \ot q_i))$ for all $x
\in A$. Then we have for every $x \in A$ that
$$\sum_{i=1}^n \eta_i(x r_i) = \sum_{i=1}^n (\th_1 \od \th_2)(\de(x)
\de(r_i)(p_i \ot q_i)) = (\th_1 \od \th_2)(\de(x)(a_1 \ot a_2))
= ((a_1 \th_1) (a_2 \th_2))(x) \ $$
which implies that $\sum_{i=1}^n r_i \, \eta_i = (a_1 \th_1) (a_2
\th_2)$. Hence, we get by definition for every $b \in B$ that
\begin{eqnarray*}
& & ((a_1 \th_1)(a_2 \th_2) \od \io)(V) \, b
= \sum_{i=1}^n (\eta_i \od \io)(V(r_i \ot b))
= \sum_{i=1}^n (\th_1 \od \th_2 \od \io)( (\de \od \io)(V(r_i \ot b))
(p_i \ot q_i \ot 1)) \\
& & \spat = \sum_{i=1}^n (\th_1 \od \th_2 \od \io)( (\de \od \io)(V)
(\de(r_i)(p_i \ot q_i) \ot b))
= (\th_1 \od \th_2 \od \io)( (\de \od \io)(V) (a_1 \ot a_2 \ot b)) \\
& & \spat = (a_1 \th_1 \od a_2 \th_2 \od \io)((\de \od \io)(V)) \, b \ .
\end{eqnarray*}
So we see that $((a_1 \th_1)(a_2 \th_2) \od \io)(V) = (a_1 \th_1 \od a_2
\th_2 \od \io)((\de \od \io)(V))$. The result follows by linearity.
\end{demo}

\begin{lemma} \label{lem4.2}
Consider a non-degenerate algebra $B$ and $V,W \in M(A \od B)$. Then we
have for every $\om_1,\om_2 \in A^\circ$ that
$(\om_1 \od \io)(V)\,(\om_2 \od \io)(W) = (\om_1 \od \om_2 \od
\io)(V_{13} \, W_{23})$.
\end{lemma}
\begin{demo}
Choose $a_1,a_2 \in A$ and $\th_1,\th_2 \in A'$. Take $b \in B$. Then we
get by definition that
\begin{eqnarray*}
(a_1 \th_1 \od \io)(V)\,(a_2 \th_2 \od \io)(W)\,b & = & (a_1 \th_1 \od
\io)(V) \, (\th_2 \od \io)(W(a_2 \ot b)) \\ & = & (\th_1 \od
\io)\bigl(\,V(a_1 \ot (\th_2 \od \io)(W(a_2 \ot b)))\,\bigr)
\hspace{1.5cm} \text{(*)}\\
\end{eqnarray*}
We have for every $p,q \in A$ that
\begin{eqnarray*}
V (a_1 \ot (\th_2 \od \io)(p \ot q)) & = &
\th_2(p) \, V(a_1 \ot q) = (\io \od \th_2 \od \io)([V(a_1 \ot q)]_{13}
(1 \ot p \ot 1)) \\
& = & (\io \od \th_2 \od \io)(V_{13}(a_1 \ot p \ot q)) \ .
\end{eqnarray*}
Using this in equality (*), we get that
\begin{eqnarray*}
& & (a_1 \th_1 \od \io)(V) \, (a_2 \th_2 \od \io)(W) \, b
=  (\th_1 \od \io)(\,(\io \od \th_2 \od \io)(V_{13}(a_1 \ot
W(a_2 \ot b)))\,) \\
& & \spat =(\th_1 \od \th_2 \od \io)(V_{13} W_{23} (a_1 \ot a_2 \ot b))
= (a_1 \th_1 \od a_2 \th_2 \od \io)(V_{13} W_{23}) \, b \ .
\end{eqnarray*}
Hence, we see that $(a_1 \th_1 \od \io)(V) \, (a_2 \th_2 \od \io)(W)
= (a_1 \th_1 \od a_2 \th_2 \od \io)(V_{13} W_{23})$.

The lemma follows.
\end{demo}

\begin{definition}
Consider a non-degenerate algebra $B$ and a corepresentation $V$ of
$(A,\de)$ on $B$. Then we define the mapping $\pi_V$ from $\ah$ into
$M(B)$ such that $\pi_V(\om) = (\om \od \io)(V)$ for every $\om \in \ah$.
Then $\pi_V$ is an algebra homomorphism.
\end{definition}

\medskip

\begin{remark} \rm \label{rem4.1}
\begin{itemize}
\item The multiplicativity follows immediately from lemmas \ref{lem4.1}
and \ref{lem4.2}. Combining these lemmas with result \ref{res3.3}, we
have also a converse:

\smallskip\smallskip

Consider a non-degenerate algebra $B$ and an element $V \in M(A \od B)$.
Then $V$ is a corepresentation $\Leftrightarrow$ We have for every
$\om_1, \om_2 \in \ah$ that $(\om_1 \od \io)(V) \, (\om_2 \od \io)(V) =
(\om_1 \om_2 \od \io)(V)$.
\item Let $B$ be a non-degenerate algebra and $V,W$ corepresentations of
$(A,\de)$ on $B$. Then result \ref{res3.3} implies immediately that $V=W$
if and only if $\pi_V=\pi_W$.
\end{itemize}
\end{remark}

\begin{proposition} \label{prop4.4}
Consider a non-degenerate algebra $B$ and a  corepresentation $V$ of
$(A,\de)$ on $B$ such that $A \od B = V (A \od B) = (A \od B) V$.
Then $\pi_V$ is non-degenerate.
\end{proposition}
\begin{demo}
Choose $b \in B$. There exists $a \in A$ such that $\vfi(a) = 1$. By
assumption, there exist $p_1,\ldots\!,p_n \in A$ and  $q_1,\ldots\!,q_n
\in B$ such that $a \ot b = \sum_{i=1}^n V(p_i \ot q_i)$. This implies that
\begin{eqnarray*}
\sum_{i=1}^n \pi_V(p_i \vfi) \, q_i & = & \sum_{i=1}^n (p_i \vfi \od
\io)(V) \, q_i =  \sum_{i=1}^n (\vfi \od \io)(V(p_i \ot q_i))  \\
& = & (\vfi \od \io)(a \ot b) = b \ .
\end{eqnarray*}

\medskip

Hence, we see that $\pi_V(\ah) B = B$. Similarly, we get that $B \,
\pi_V(\ah) = A$.
\end{demo}

Notice that the conclusion of this proposition is especially true for a
non-degenerate corepresentation. Later on, we will prove that the
non-degeneracy of $V$ is equivalent with the non-degeneracy of $\pi_V$.

\bigskip

Next we prove some results about the behaviour of a corepresentation with
respect to the counit and the antipode.

\medskip

\begin{proposition}
Consider a non-degenerate algebra $B$ and a non-degenerate
corepresentation $V$ of $(A,\de)$ on $B$. Then $(\vep \od \io)(V) = (\io
\od \vep)(V) = 1$.
\end{proposition}
\begin{demo}
We have that $(\de \od \io)(V) = V_{13} V_{23}$. Applying $\vep \od \io
\od \io$ to this equation implies that
$$(\,(\vep \od \io)\de \od \io)(V) = (1 \ot (\vep \od \io)(V)) \,\, V \ .$$
Using the fact that $(\vep \od \io)\de = \io$, this implies that
$V = (1 \ot (\vep \od \io)(V)) \, \, V$. Hence, the invertibility of $V$
implies that $1 \ot (\vep \od \io)(V) = 1 \ot 1$, so $(\vep \od
\io)(V)=1$. Similarly, one proves that $(\io \od \vep)(V)=1$.
\end{demo}

\medskip

\begin{proposition} \label{prop4.1}
Consider a non-degenerate algebra $B$ and a non-degenerate
corepresentation $V$ of $(A,\de)$ on $B$. Then we have the following
equalities:
\begin{enumerate}
\item $A \od B = \langle \, (a \ot 1) V (1 \ot b) \mid a \in A, b \in B
\,\rangle$
\item $A \od B = \langle \, (a \ot 1) V^{-1} (1 \ot b) \mid a \in A, b
\in B \,\rangle$
\item $A \od B = \langle \, (1 \ot b) V (a \ot 1) \mid a \in A, b \in B
\,\rangle$
\item $A \od B = \langle \, (1 \ot b) V^{-1} (a \ot 1) \mid a \in A, b
\in B \,\rangle$
\end{enumerate}
\end{proposition}
\begin{demo}
Choose $a,c \in A$ and $b \in B$. We are going to prove the equality
$$(\psi \od \io \od \io)((V^{-1})_{13}\,((1 \ot a)\de(c) \ot b))
= (a \ot 1) V (1 \ot (\psi \od \io)(V^{-1}(c \ot b))\,) \ . $$

\medskip

Choose $d \in B$. Then we have that
\begin{eqnarray*}
& & \hspace{-0.5cm} (1 \ot d)\,(\psi \od \io \od \io)((V^{-1})_{13}
\,((1 \ot a)\de(c) \ot b))
= (\psi \od \io \od \io)(\,(1 \ot 1 \ot d) (V^{-1})_{13} (1 \ot a \ot
1)(\de(c) \ot b)) \\
& & \hspace{-0.78cm} \spat = (\psi \od \io \od \io)( (1 \ot a \ot d)
(V^{-1})_{13} (\de(c) \ot b))
= (\psi \od \io \od \io)( (1 \ot a \ot d) V_{23} (\de \od \io)(V^{-1})
(\de \od \io)(c \ot b))
\end{eqnarray*}
where in the last equality, we used the fact that $V$ is a
corepresentation. So we get that
$$ (1 \ot d)\,(\psi \od \io \od \io)((V^{-1})_{13}((1 \ot a)\de(c) \ot
b)) = (\psi \od \io \od \io)(\,(1 \ot (a \ot d)V)(\de \od \io)(V^{-1}(c
\ot b))\,) \ .  \text{\ \ \ \ (*)} $$

\begin{list}{}{\setlength{\leftmargin}{.4 cm}}

\item We have for every $p,r \in A$ and $q,s \in B$ that
\begin{eqnarray*}
& & (\psi \od \io \od \io)(\,(1 \ot p \ot q) (\de \od \io)(r \ot s)\,)
= (\psi \od \io \od \io)(\,(1 \ot p)\de(r) \ot q s\,) \\
& & \spat = \psi(r) \,p \ot q s = (p \ot q)(1 \ot (\psi \od \io)(r \ot
s)) .
\end{eqnarray*}
where we used the right invariance of $\psi$ in the second last equality.

\end{list}

Combining this with equality (*), we arrive at the conclusion that
$$(1 \ot d)\,(\psi \od \io \od \io)((V^{-1})_{13}((1 \ot a)\de(c) \ot b))
= (a \ot d) V \, (1 \ot (\psi \od \io)(V^{-1}(c \ot b))\,) .$$
So we get that
$$(\psi \od \io \od \io)((V^{-1})_{13}((1 \ot a)\de(c) \ot b)) = (a \ot
1) V \, (1 \ot (\psi \od \io)(V^{-1}(c \ot b))\,) \ .$$

Using this equality, it is not difficult to check that $A \od B =
\langle \, (a \ot 1) V (1 \ot b) \mid a \in A, b \in B \, \rangle$. The
other equalities are proven in a similar way (The last two, by using
the left Haar functional).
\end{demo}

Now we get easily the following result:

\begin{result}
Consider a non-degenerate algebra $B$ and a non-degenerate
corepresentation $V$ of $(A,\de)$ on $B$. Let $a$ be an element in $A$
and $\om$ be an element in $A'$. Then $(a \om  \od \io)(V)$ belongs to
$M(B)$ and $b \, (a \om \od \io)(V) = (\om \od \io)( (
1 \ot b) V (a \ot 1))$ for every $b \in B$.
\end{result}

Proposition \ref{prop4.1} makes it also possible to give the following
definition.

\begin{definition}
Consider a non-degenerate algebra $B$ and a non-degenerate
corepresentation $V$ of $(A,\de)$ on $B$. Then $(S \od \io)(V)$ is
defined to be the unique element in $M(A \od B)$ such that \newline
$(S \od \io)(V) \, (a \ot b) = (S \od \io)((S^{-1}(a) \ot 1) V (1 \ot b))$
and
$(a \ot b) \, (S \od \io)(V) = (S \od \io)( (1 \ot b) V (S^{-1}(a)
\ot 1))$ for every $a \in A$ and $b \in B$.
\end{definition}

It is not difficult to check that $(S \od \io)(V)$ is really an
element of $M(A \ot B)$.

\medskip

In this terminology, we get the usual result:

\begin{proposition} \label{prop4.3}
Consider a non-degenerate algebra $B$ and a non-degenerate
corepresentation $V$ of $(A,\de)$ on $B$. Then $(S \od \io)(V) =
V^{-1}$.
\end{proposition}
\begin{demo}
Choose $b \in B$ and $c,d \in A$. Take $a \in A$ such that $\vep(a)=1$.

Because $(\de \od \io)(V) = V_{13} V_{23}$, we have the equality
\begin{eqnarray*}
(\de \od \io)((a \ot 1)V(1 \ot b)) \, (c \ot d \ot 1)
& = & (\de(a) \ot 1) V_{13} V_{23} (c \ot d \ot b) \\
& = & (\de(a) \ot 1) [V(c \ot 1)]_{13} [V(d \ot b)]_{23} \hspace{1.5cm}
(a) \\
\end{eqnarray*}
By  proposition \ref{prop4.1}, we know that $(a \ot 1)V(1 \ot b)$
belongs to $A \od B$. It is not difficult to see that $[V(c \ot
1)]_{13} [V(d \ot b)]_{23}$ belongs to $A \od A \od B$.

\medskip

It is rather easy to check that $m(S \od \io)(\de(x) \, y) = \vep(x)
\, m(S \od \io)(y)$ for every $x,y \in A \od A$.

Therefore, if we apply $m(S \od \io) \, \od \io$ to equation (a), we
get the equality
$$m(S \od \io)(c \ot d) \ot (\vep \od \io)((a \ot 1)V(1 \ot b))
= \vep(a) \, (m(S \od \io) \, \od \io)([V(c \ot 1)]_{13} [V(d \ot
b)]_{23}) \ .\hspace{1cm} (b)$$
Using the fact that $(\vep \od \io)(V)=1$, the left hand side of equation
(b) is equal to $S(c) d \ot b$.

\medskip

Now we are going to work on the right hand side.
We have for every $p,q \in A$ that
\begin{eqnarray*}
(m(S \od \io) \, \od \io)([V(c \ot 1)]_{13} [p \ot q]_{23}) & = &
(m(S \od \io) \, \od \io)([V(c \ot q)]_{13} (1 \ot p \ot 1)) \\
& = & (S \od \io)(V(c \ot q)) \, (p \ot 1)
= (S(c) \ot 1) \, (S \od \io)(V) \, (p \ot q) \ . \\
\end{eqnarray*}
This implies that the right hand side of (b) is equal to
$(S(c) \ot 1) \, (S \od \io)(V) \, V(d \ot b)$.

\medskip

So we get the equality $S(c) d \ot b = (S(c) \ot 1) \, (S \od \io)(V) \,
V(d \ot b)$.

\medskip

Consequently, we have proven that $(S \od \io)(V) \, V = 1$ which implies
that $(S \od \io)(V) = V^{-1}$.
\end{demo}

\bigskip

If $\pi_V$ is non-degenerate, we can extend it to $M(\ah)$ which
contains $A^\circ$. For these elements, $\pi_V$ will act in the
obvious way:

\begin{proposition}
Consider a non-degenerate algebra $B$ and a non-degenerate
corepresentation $V$ of  $(A,\de)$ on $B$. Then we have for every $a \in
A$, $\om \in A'$ and $b \in B$ that
$b \, \pi_V(a \om)  = (\om \od \io)((1 \ot b)V(a \ot 1))$ and
$\pi_V(a \om) \, b = (\om \od \io)(V (a \ot b))$.
\end{proposition}
\begin{demo}
Remember that $a \om$ is an element of $A^\circ$, so we have already
defined the element $(a \om \od \io)(V)$. Choose $\th
\in \ah$ and $b \in B$. Using lemmas \ref{lem4.1} and \ref{lem4.2}, we get
that
\begin{eqnarray*}
& & \pi_V(a \om)(\pi_V(\th) b)  = \pi_V((a \om) \th) \, b
= ((a \om) \th \od \io)(V) \, b =  ((a \om) \od \th \od \io)((\de \od
\io)(V)) \, b \\
& & \spat = (a \om \od \th \od \io)(V_{13} V_{23})\, b
=  (a \om \od \io)(V) \, (\th \od \io)(V) \,b
= (a \om \od \io)(V) \, \pi_V(\th) b \ .
\end{eqnarray*}
Hence, the non-degeneracy of $\pi_V$ implies that $\pi_V(a \om) = (a
\om \od \io)(V)$.
\end{demo}

Of course, a similar result is true for linear functionals  of the
form $\om a$. In this case, the results in the previous section would
have to be stated in terms of right multipliers.

\medskip\medskip

\begin{proposition}  \label{prop4.2}
Consider non-degenerate algebras $B,C$ and a non-degenerate homomorphism
$\th$ from $B$ into $M(C)$. Let $V$ be a corepresentation of $(A,\de)$ on
$B$. Then $(\io \od \th)(V)$ is a corepresentation of $(A,\de)$ on $C$
such that $\pi_{(\io \od \th)(V)} =  
\overline{\th} \circ \pi_V$. If $V$ is non-degenerate, then $(\io \od
\th)(V)$ is non-degenerate.
\end{proposition}
\begin{demo}
We have that
\begin{eqnarray*}
& & (\de \od \io)((\io \od \th)(V)) = (\io \od \io \od \th)((\de \od \io)(V))
= (\io \od \io \od \th)(V_{13} V_{23}) \\
& & \spat = (\io \od \io \od \th)(V_{13}) \, (\io \od \io \od \th)(V_{23})
= (\io \od \th)(V)_{13} \, (\io \od \th)(V)_{23}
\end{eqnarray*}
which implies that $(\io \od \th)(V)$ is a corepresentation.

Using \ref{res3.2}, we have for every $\om \in \ah$ that
$$\pi_{(\io \od \th)(V)}(\om) = (\om \od \io)((\io \od \th)(V))
= \th( (\om \od \io)(V)) = \th(\pi_V(\om)) \ . $$
\end{demo}

In  the following section, we show the existence of a special
corepresentation such that every corepresentation can be constructed
from this special one using the method described in the proposition
above.

\section{The universal corepresentation of an algebraic quantum group}
\label{art5}

Consider an algebraic quantum group $(A,\de)$ with a left Haar functional
$\vfi$. We will construct the universal corepresentation of $(A,\de)$ and
use it to prove that there is a natural bijection between
corepresentations of $(A,\de)$ and representations of $\ah$.

\medskip

Because $\ah$ is a subset of $A'$, we can regard $A \od \ah$ in a natural
way as a subspace of the linear operators on $A$.

\begin{notation} \label{not5.1}
We define linear mappings $U_l,U_r$ from $A \od \ah$ into $L(A)$ such
that  $[U_l(x \ot \om)](y) = (\io \od \om)(\de(y)(x \ot 1))$
and $[U_r(x \ot \om)](y) = (\om \od \io)((1 \ot x)\de(y))$
for every $x,y \in A$ and $\om \in \ah$
\end{notation}

\begin{proposition} \label{prop5.1}
We have for every $a \in A \od \ah$ that $U_l(a)$ and $U_r(a)$ belong to
$A \od \ah$.
\end{proposition}
\begin{demo}
Choose $x \in A$ and $\om \in \ah$. Then there exists $z \in A$ such that
$\om = z \vfi$. Furthermore, there exist $p_1,\ldots\!,p_n ,
q_1,\ldots\!,q_n \in A$ such that $x \ot z = \sum_{i=1}^n \de(p_i)(q_i \ot 1)$.

Then we have for every $y \in A$ that
\begin{eqnarray*}
[U_l(x \ot \om)] (y) & = & (\io \od \om)(\de(y)(x \ot 1))
= (\io \od \vfi)(\de(y)(x \ot z)) \\
& = & \sum_{i=1}^n (\io \od \vfi)(\de(y p_i) (q_i \ot 1))
= \sum_{i=1}^n \vfi(y p_i) \, q_i
\end{eqnarray*}
Where in the last equality, we used the left invariance of $\vfi$.
So we get that $U_l(x \ot \om) = \sum_{i=1}^n q_i \ot p_i \vfi$.

Similarly, we get that $U_r(x \ot \om)$ belongs to $A \od \ah$.
\end{demo}

As usual, $A \od \ah$ has a natural algebra structure. In a next step, we
want to prove that $(U_l,U_r)$ is a two-sided multiplier on $A \od \ah$.
First we need an easy lemma.

\begin{lemma}
Consider $\om \in A \od \ah$, $a,b \in A$ and $x \in A$. Then we have
for every $x \in A$ that  \newline $((b \ot \vfi a) \om)(x) = m( (b
\ot \vfi) \od \om)((a \ot 1)\de(x))$ and $(\om (b \ot \vfi a))(x) = m(\om
\od (b \ot \vfi))( (1 \ot a)\de(x))$.
\end{lemma}
\begin{demo}
Choose $c \in A$ and $\th \in \ah$. Take $x \in A$, then we have that
$$[(b \ot \vfi a)(c \ot \th)](x) = [(b c) \ot ((\vfi a) \th)](x)
= b c \,\, ((\vfi a)\th)(x) = b c \,\, (\vfi \od \th)((a \ot 1)\de(x))  \ .$$
We have for every $p,q \in A$ that
$$b c \,\, (\vfi \od \th)(p \ot q) =  b \, \vfi(p) \,\, c \, \th(q)
= (b \ot \vfi)(p) \, (c \ot \th)(q) = m( (b \ot \vfi)
\od (c \ot \th))(p \ot q) \ . $$
This implies that
$$[(b \ot \vfi a)(c \ot \th)](x) = m( (b \ot \vfi) \od (c \ot \th) ) ((a
\ot 1)\de(x)) \ .$$
The result follows by linearity. The other equality is proven in the
same way.
\end{demo}

\begin{lemma}
Consider $\om,\th \in A \od \ah$. Then
\begin{enumerate}
\item $U_l(\om) \, \th = U_l(\om \th)$
\item $\om \, U_r(\th) = U_r(\om \th)$
\item $U_l(\om) \, \th = \om \, U_r(\th)$
\end{enumerate}
\end{lemma}
\begin{demo}
Choose $a_1, a_2, b_1, b_2 \in A$.
\begin{enumerate}
\item Take $x \in A$. Then we have by the previous lemma that
$$[U_l(b_1 \ot \vfi a_1) \, (b_2 \ot \vfi a_2)](x)
= m( U_l(b_1 \ot \vfi a_1) \ot (b_2 \ot \vfi))((1 \ot a_2)\de(x)) \ .
\hspace{1.5cm} \text{(a)}$$

\medskip

We have for every $p,q \in A$ that
\begin{eqnarray*}
& & m( U_l(b_1 \ot \vfi a_1) \ot (b_2 \ot \vfi))(p \ot q)
 = [U_l(b_1 \ot \vfi a_1)](p) \, (b_2 \ot \vfi)(q) \\
& & \spat = (\io \od \vfi a_1)(\de(p)(b_1 \ot 1)) \, b_2 \vfi(q)
= (\io \od \vfi a_1 \od \vfi)( \de(p)(b_1 b_2 \ot 1) \ot q) \\
& & \spat = (\io \od \vfi a_1 \od \vfi)((\de \od \io)(p \ot q)(b_1 b_2
\ot 1 \ot 1))  \ .
\end{eqnarray*}

\medskip

Using this with equation (a), we get that
\begin{eqnarray*}
& & [U_l(b_1 \ot \vfi a_1) \, (b_2 \ot \vfi a_2)](x)
= (\io \od \vfi a_1 \od \vfi)((\de \od \io)((1 \ot a_2)\de(x))(b_1 b_2
\ot 1 \ot 1))  \\
& & \spat = (\io \od \vfi a_1 \od \vfi)( (1 \ot 1 \ot a_2) (\de \od
\io)(\de(x)) (b_1 b_2 \ot 1 \ot 1) ) \\
& & \spat = (\io \od \vfi a_1 \od \vfi)( (1 \ot 1 \ot a_2) (\io \od
\de)(\de(x)) (b_1 b_2 \ot 1 \ot 1)) \\
& & \spat = (\io \od \vfi a_1 \od \vfi)( (1 \ot 1 \ot a_2) (\io \od
\de)(\de(x)(b_1 b_2 \ot 1)) \, )  \ .
\end{eqnarray*}
Using definition \ref{def1.1}, this implies that
\begin{eqnarray*}
& & [U_l(b_1 \ot \vfi a_1) \, (b_2 \ot \vfi a_2)](x)
= (\io \od (\vfi a_1)(\vfi a_2))(\de(x)(b_1 b_2 \ot 1))  \\
& & \spat = [U_l( b_1 b_2 \ot (\vfi a_1)(\vfi a_2) )](x)
= [U_l( (b_1 \ot \vfi a_1)(b_2 \ot \vfi a_2))](x) \ .
\end{eqnarray*}
So we get that $U_l(b_1 \ot \vfi a_1) \, (b_2 \ot \vfi a_2)
=  U_l( (b_1 \ot \vfi a_1)(b_2 \ot \vfi a_2))$.
\item In a similar way, one proves that
$(b_1 \ot \vfi a_1) \, U_r(b_2 \ot \vfi a_2)
=  U_r( (b_1 \ot \vfi a_1)(b_2 \ot \vfi a_2))$.
\item Take $x \in A$. Again, the previous lemma implies that
$$ [U_r(b_1 \ot \vfi a_1)\,(b_2 \ot \vfi a_2)](x)
= m( U_r(b_1 \ot \vfi a_1) \ot (b_2 \ot \vfi))( (1 \ot a_2)\de(x) ) \ .
\hspace{1.5cm} \text{(b)} $$

\medskip

We have for every $p,q \in A$ that
\begin{eqnarray*}
& & m( U_r(b_1 \ot \vfi a_1) \ot (b_2 \ot \vfi))(p \ot q)
= [U_r(b_1 \ot \vfi a_1)](p) \, (b_2 \ot \vfi)(q) \\
& & \spat = (\vfi a_1 \ot \io)((1 \ot b_1)\de(p)) \, b_2 \vfi(q)
= (\vfi a_1 \ot \io)((1 \ot b_1)\de(p)(1 \ot b_2)) \, \vfi(q) \\
& & \spat = (\vfi a_1 \od \io \od \vfi)( (1 \ot b_1)\de(p)(1 \ot b_2) \ot
q )\\
& & \spat = (\vfi a_1 \od \io \od \vfi)( (1 \ot b_1 \ot 1)(\de \od \io)(p
\ot q)(1 \ot b_2 \ot 1))  \ .
\end{eqnarray*}

\medskip

Using this, equation (b) implies that
$$[U_r(b_1 \ot \vfi a_1)\,(b_2 \ot \vfi a_2)](x)
= (\vfi a_1 \od \io \od \vfi)( (1 \ot b_1 \ot 1)(\de \od \io)((1 \ot
a_2)\de(x))(1 \ot b_2 \ot 1))  \ , $$
which implies that
$$[U_r(b_1 \ot \vfi a_1)\,(b_2 \ot \vfi a_2)](x)
= (\vfi \od \io \od \vfi)( (a_1 \ot b_1 \ot a_2) (\de \od \io)(\de(x))
(1 \ot b_2 \ot 1)) \ .  \hspace{1cm} \text{(c)}$$

In a similar way, we get that
$$[(b_1 \ot \vfi a_1)\,U_l(b_2 \ot \vfi a_2)](x)
= (\vfi \od \io \od \vfi a_2)( (1 \ot b_1 \ot 1) (\io \od \de)((a_1 \ot
1)\de(x)) (1 \ot b_2 \ot 1)) \ , $$
which implies that
$$[(b_1 \ot \vfi a_1)\,U_l(b_2 \ot \vfi a_2)](x)
= (\vfi \od \io \od \vfi)( (a_1 \ot b_1 \ot a_2) (\io \od \de)(\de(x))
(1 \ot b_2 \ot 1)) \ . \hspace{1cm} \text{(d)}$$

Hence, using the coassociativity of $\de$ and looking at equalities (c)
and (d), we arrive at the conclusion that
$$[U_r(b_1 \ot \vfi a_1)\,(b_2 \ot \vfi a_2)](x)
=[(b_1 \ot \vfi a_1)\,U_l(b_2 \ot \vfi a_2)](x) \ .$$

So we have proven that $U_r(b_1 \ot \vfi a_1)\,(b_2 \ot \vfi a_2)
=(b_1 \ot \vfi a_1)\,U_l(b_2 \ot \vfi a_2)$.
\end{enumerate}
\end{demo}

This lemma justifies the following definition:

\begin{definition} \label{def5.1}
There exists a unique element $U \in M(A \od \ah)$ such that $[U(x \ot
\om)](y) =$ \newline $(\io \od \om)(\de(y)(x \ot 1))$ and $[(x \ot
\om)U](y) = (\om \od \io)((1 \ot x)\de(y))$ for every $x,y \in A$ and
$\om \in \ah$.
\end{definition}

Later on, we will prove that $U$ is a corepresentation of $(A,\de)$ on
$\ah$. It is called the universal corepresentation of $(A,\de)$.

\medskip

\begin{proposition}
The element $U$ is invertible in $M(A \od \ah)$.
\end{proposition}
\begin{demo}
We prove that $U_l$ is bijective:
\begin{itemize}
\item First we prove that $U_l$ is injective.

Choose $x_1,\ldots\!,x_n \in A$ and $\om_1,\ldots\!,\om_n \in \ah$
such that $U_l(\sum_{i=1}^n x_i \ot \om_i) = 0$.

Take $y \in A$. Choose $b \in A$. Then there exist $p_1,\ldots\!,p_m ,
q_1,\ldots\!,q_m \in A$ such that $b \ot y = \sum_{j=1}^m (p_j \ot
1)\de(q_j)$. So we get that
\begin{eqnarray*}
0 & = & \sum_{j=1}^m p_j \, [U_l(\sum_{i=1}^n x_i \ot \om_i)](q_j)
= \sum_{j=1}^m \sum_{i=1}^n p_j \, (\io \od \om_i)(\de(q_j)(x_i \ot 1)) \\
& = & \sum_{i=1}^n \sum_{j=1}^m (\io \od \om_i)([(p_j \ot 1)\de(q_j)](x_i
\ot 1))
= \sum_{i=1}^n (\io \od \om_i)((b \ot y)(x_i \ot 1)) \\
& = & b x_i \, \om_i(y) = b \,\, (\sum_{i=1}^n x_i \ot \om_i)(y) \ .
\end{eqnarray*}
So we get that $(\sum_{i=1}^n x_i \ot \om_i)(y) = 0$. Consequenctly
$\sum_{i=1}^n x_i \ot \om_i = 0$.

\item Next we prove that $U_l$ is surjective.

Choose $a,b \in A$. Then there exist $p_1,\ldots\!,p_n ,
q_1,\ldots\!,q_n \in A$ such that $\de(b)(a \ot 1) = \sum_{i=1}^n p_i
\ot q_i$. Then we have for every $x \in A$ that
\begin{eqnarray*}
[U_l(\sum_{i=1}^n p_i \ot q_i \vfi)](x) & = &
\sum_{i=1}^n (\io \od q_i \vfi)(\de(x)(p_i \ot 1))
= \sum_{i=1}^n (\io \od \vfi)(\de(x)(p_i \ot q_i)) \\
& = & (\io \od \vfi)(\de(x b)(a \ot 1))
= \vfi(x b) a = (a \ot b \vfi)(x) \ .
\end{eqnarray*}
Consequently $a \ot b \vfi = U_l(\sum_{i=1}^n p_i \ot q_i \vfi)$.
\end{itemize}

In a similar way, one proves that $U_r$ is bijective. These two facts
imply that $U$ is invertible in $M(A \od \ah)$.
\end{demo}

\begin{lemma}
Consider $x \in A \od \ah$ and $\om \in A'$. Then
$[(\om \od \io)(x)](a) = \om(x(a))$ for every $a \in A$.
\end{lemma}
\begin{demo}
Choose $b \in A$ and $\th \in \ah$. Then we have for every $a \in A$
that
$$[(\om \od \io)(b \ot \th)](a) = [\om(b) \th](a)
= \om(b) \th(a) = \om(b \th(a)) = \om( (b \ot \th)(a)) \ . $$
The lemma follows.
\end{demo}

\begin{result}
We have for every $\om \in \ah$ that
$(\om \od \io)(U) = \om$.
\end{result}
\begin{demo}
There exists $a \in A$ such that $\om = \vfi a$.

Take $\th \in \ah$. Then we have for every $x \in A$ that
\begin{eqnarray*}
(\th \, (\vfi a \od \io)(U))(x) & = &
(\vfi \od \io)( (a \ot \th) U) (x)
\stackrel{(*)}{=} \vfi( [(a \ot \th)U](x) ) \\
& = & \vfi(\,(\th \od \io)((1 \ot a)\de(x))\, )
=  (\th (\vfi a))(x)
\end{eqnarray*}
where the previous lemma was used in the equality (*).
So we get that $\th \, (\vfi a \od \io)(U)
= \th (\vfi a)$. Therefore, lemma \ref{lem5} implies that $(\vfi a \od
\io)(U) = \vfi a$.
\end{demo}

\begin{corollary} \label{cor5.1}
The element $U$ is a non-degenerate corepresentation of $(A,\de)$ on
$\ah$ such that $\pi_U = \io$.
\end{corollary}

This follows easily from remark \ref{rem4.1} and the previous result.

\medskip

\begin{proposition}
Consider a non-degenerate algebra $B$ and a corepresentation $V$ of
$(A,\de)$ on $B$. Then $V$ is non-degenerate $\Leftrightarrow$  $\pi_V$
is non-degenerate$\Leftrightarrow$  $V(A \od B) = (A \od B) V = A \od B$
\end{proposition}
\begin{demo}
\begin{itemize}
\item If $V$ is non-degenerate, we get easily that $V(A \od B) = (A \od
B) V = A \od B$.
\item If $V(A \od B) = (A \od B) V = A \od B$, we get by \ref{prop4.4}
that $\pi_V$ is non-degenerate.
\item Suppose that $\pi_V$ is non-degenerate. By result \ref{res3.2}, we
have for every $\om \in \ah$ that
$$(\om \ot \io)( (\io \ot \pi_V)(U) ) = \pi_V( (\om \ot \io)(U) )
= \pi_V(\om) = (\om \ot \io)(V) \ ,$$
which, by result \ref{res3.3}, implies that $(\io \ot \pi_V)(U) = V$.

Because $U$ is invertible and $\pi_V$ is non-degenerate, this implies
that $V$ is invertible.
\end{itemize}
\end{demo}

In fact, we have even proven the following proposition:

\begin{proposition}
Consider a non-degenerate algebra $B$ and a non-degenerate
corepresentation $V$ of $(A,\de)$ on $B$. Then $(\io \od \pi_V)(U) = V$.
\end{proposition}

Hence, $U$ satisfies the following universal property (which was
introduced for compact quantum groups in \cite{PW}).

\begin{proposition}
Consider a non-degenerate algebra $B$ and a non-degenerate
corepresentation $V$ of $(A,\de)$ on $B$.
Then there exists a unique non-degenerate homomorphism $\th$ from $\ah$
into $M(B)$ such that $(\io \od \th)(U) = V$.
\end{proposition}
\begin{demo}
We have already proven the existence. We prove now quickly the uniqueness.

Therefore, let $\th$ be a non-degenerate homomorphism from $\ah$ into
$M(B)$ such that $(\io \od \th)(U) = V$.

Then we have for every $\om \in \ah$ that
$$\pi_V(\om) = (\om \od \io)(V) = (\om \od \io)((\io \ot \th)(U))
= \th((\om \ot \io)(U)) = \th(\om) \, $$
So we see that $\pi_V = \th$.
\end{demo}

\medskip

Using proposition \ref{prop4.2}, we have also a converse of the previuous
result.

\begin{proposition}
Consider a non-degenerate algebra $B$ and a non-degenerate homomorphism
$\th$ from $\ah$ into $M(B)$. Then $(\io \od \th)(U)$ is a non-degenerate
corepresentation of $(A,\de)$ on $B$ and
$\pi_{(\io \od \th)(U)} = \th$.
\end{proposition}

\medskip

Consequently, we have proven that there is a bijective correspondence
between non-degenerate corepresentations of $(A,\de)$ and homomorphisms
on $\ah$. The universal corepresentation serves as a linking mechanism.

\medskip

\begin{result}
Consider a non-degenerate algebra $B$ and a non-degenerate
corepresentation $V$ of $(A,\de)$ on $B$. Then
$(\om \od \io)(V^{-1}) = \pi_V(\hat{S}(\om))$ for every $\om \in \ah$.
\end{result}
\begin{demo}
Choose $a \in A$. Using proposition \ref{prop4.3}, we have for every $b
\in B$ that
\begin{eqnarray*}
& & (a \vfi  \od \io)(V^{-1}) \, b
= (\vfi \od \io)(V^{-1} (a \ot b)) = (\vfi \od \io)(\,(S \od \io)(
(S^{-1}(a) \ot 1) V (1 \ot b))\,) \\
& & \spat = (\psi \od \io)( (S^{-1}(a) \ot 1) V (1 \ot b) )
= (\psi \,S^{-1}(a) \od \io)(V) \, b
= \pi_V(\psi\,S^{-1}(a)) \, b \ .
\end{eqnarray*}
Because $\psi\,S^{-1}(a) = (a \vfi)\!\circ\!S = \hat{S}(a \vfi)$, we get
that
$$(a \vfi  \od \io)(V^{-1}) = \pi_V(\psi\,S^{-1}(a)) = \pi_V(\hat{S}(a
\vfi)) \ .$$
\end{demo}

\section{Unitary corepresentations of $^*$-algebraic quantum groups}

We start this section with the definition of a $^*$-algebraic quantum group.

\begin{definition}
Consider a non-degenerate $^*$-algebra $A$ and a non-degenerate
$^*$-homomorphism $\de$
from $A$ into $M(A \od A)$ such that
\begin{enumerate}
\item $(\de \od \io)\de = (\io \od \de)\de$.
\item The linear mappings $T_1$, $T_2$ from $A \od A$ into $M(A \od A)$
such that
$$T_1(a \ot b) = \de(a)(b \ot 1) \text{\ \ \ \ and \ \ \ \ } T_2(a
\ot b) = \de(a)(1 \ot       b)$$
for all $a,b \in A$, are bijections from $A \od A$ to $A \od A$.
\end{enumerate}
Then we call $(A,\de)$ a Multiplier Hopf$\,^*$-algebra.
\end{definition}

It is easy to see that a Multiplier Hopf$\,^*$-algebra is a regular
Multiplier Hopf algebra.

\medskip

\begin{definition}
Consider a Multiplier Hopf$\,^*$-algebra $(A,\de)$ such that there exists
a non-zero linear functional $\vfi$ on $A$ which is left invariant.  Then
we call $(A,\de)$ a $^*$-algebraic quantum group.
\end{definition}

So every $^*$-algebraic quantum group is an algebraic quantum group.

\medskip

\begin{remark} \rm Consider a $^*$-algebraic quantum group $(A,\de)$. Let
$\vfi$ be a left Haar functional of $(A,\de)$. Because $\de$ is a
$^*$-homomorphism, $\overline{\vfi}$ will also be left invariant.
So $\frac{\vfi + \overline{\vfi}}{2}$ and $\frac{\vfi -
\overline{\vfi}}{2i}$ are self adjoint left invariant functionals.
Because their sum is equal to $\vfi$, one of them has to be non zero.

\smallskip

Hence, we get the existence of a self adjoint left Haar functional on
$(A,\de)$.
\end{remark}

\bigskip\medskip

For the rest of this section, we fix a $^*$-algebraic quantum group
$(A,\de)$ with a self adjoint left Haar functional $\vfi$.
In this case $M(\ah)$ is also a $^*$-algebra. First we prove the usual
formula for the adjoint operation on $M(\ah)$.

\begin{proposition}
Let $\om$ be an element of $M(\ah)$. Then $\om^*(x) =
\overline{\om(S(x)^*)}$ for every $x \in A$.
\end{proposition}
\begin{demo}
It is not very difficult to check that $(\vfi b)^* = \psi S(b)^*$ and
$(\psi b)^* = \vfi S(b)^*$ for every $b \in A$.

Choose $a \in A$.
Then lemma \ref{lem2} implies that
$$ \om^* \, (\vfi a) = \vfi \, [S^{-1}\bigl((\io \od
\om^*)\de(S(a))\bigr)] $$
On the other hand, we have that
$$ \om^* \, (\vfi a) =  ( (\vfi a)^* \om )^*
= ( (\psi S(a)^*) \om )^* \, $$
By lemma \ref{lem2}, we know that $(\psi S(a)^*) \om = \psi\,[S\bigl((\om
\od \io)\de(S^2(a)^*)\bigr)]$, so
$$ \om^* (\vfi a) = ( \psi \, [S\bigl((\om \od
\io)\de(S^2(a)^*)\bigr)]\,)^*
=  \vfi \, [S^2\bigl((\om \od \io)\de(S^2(a)^*)\bigr)^*] \ .$$
Again, the faithfulness of $\vfi$ implies that
$$S^{-1}\bigl((\io \od \om^*)\de(S(a))\bigr) = S^2(\bigl(\om \od
\io)\de(S^2(a)^*)\bigr)^* . $$
Applying $\vep$ to this equation and using lemma \ref{lem3}, we get that
$$\om^*(S(a)) = \overline{\om(S^2(a)^*)} = \overline{\om(S(S(a))^*)} \ .$$
\end{demo}

\begin{definition}
Consider a non-degenerate $^*$-algebra $B$. A unitary corepresentation of
$(A,\de)$ on $B$ is by definition a corepresentation of $(A,\de)$ on $B$
which is a unitary element in the $^*$-algebra $M(A \od B)$.
\end{definition}

It is clear that a unitary corepresentation is automatically non-degenerate.

\begin{remark} \rm
Consider non-degenerate $^*$-algebras $B$, $C$ and a non-degenerate
$^*$-homomorphism $\th$ from $B$ into $M(C)$. If $V$ is a unitary
corepresentation of $(A,\de)$ on $B$, then it is clear that
$(\io \od \th)(V)$ is a unitary corepresentation of $(A,\de)$ on $C$.
\end{remark}

\medskip

\begin{proposition}
Consider a non-degenerate $^*$-algebra $B$ and a non-degenerate
corepresentation $V$ of $(A,\de)$ on $B$. Then
$(\om \od \io)(V^*) = \pi_V(\overline{\om}\,)^*$ for every $\om \in \ah$.
\end{proposition}
\begin{demo}
Choose $a \in A$. Then we have for every $b \in B$ that
\begin{eqnarray*}
& & (a \vfi \od \io)(V^*) \, b
= (\vfi \od \io)(V^* (a \ot b))
= (\vfi \od \io)((a^* \ot b^*) V)^* \\
& & \spat = [b^* \, (\vfi a^* \od \io)(V)]^*
= (\overline{a \vfi} \od \io)(V)^* \, b
= \pi_V(\overline{a \vfi})^* \, b \ .
\end{eqnarray*}
So we get that $(a \vfi \od \io)(V^*) = \pi_V(\overline{a \vfi})^*$.
\end{demo}

\begin{proposition}
Consider a non-degenerate $^*$-algebra $B$ and a non-degenerate
corepresentation $V$ of $(A,\de)$ on $B$. Then $V$ is unitary
$\Leftrightarrow$ $\pi_V$ is a $^*$-homomorphism.
\end{proposition}
\begin{demo}
\begin{itemize}
\item First suppose that $V$ is unitary. Choose $\om \in \ah$. Then
$$\pi_V(\om)^* = (\overline{\om} \od \io)(V^*)
= (\overline{\om} \od \io)(V^{-1}) = \pi_V(\hat{S}(\overline{\om}))
= \pi_V(\om^*) \ . $$
So we see that $\pi_V$ is a $^*$-homomorphism.
\item Next suppose that $\pi_V$ is a $^*$-homomorphism. Choose $\om \in
\ah$. Then
$$(\om \od \io)(V^*) = \pi_V(\overline{\om})^*
= \pi_V(\overline{\om}^*) = \pi_V(\hat{S}(\om)) = (\om \od \io)(V^{-1})
\ . $$
Therefore, result \ref{res3.3} implies that $V^* = V^{-1}$.
\end{itemize}
\end{demo}

\begin{corollary}
The universal corepresentation $U$ is a unitary corepresentation of
$(A,\de)$ on $\ah$.
\end{corollary}

This follows immediately from the fact that $\pi_U$ is the identity
mapping which is obviously a \newline $^*$-homomorphism.

\medskip

This result together with the results from the  previous section imply
that there is a bijective correspondence between unitary
corepesentations on $(A,\de)$ and non-degenerate $^*$-homomorphisms on
$\ah$.

\section{The universal corepresentation of the dual}

Consider an algebraic quantum group $(A,\de)$ with a left Haar
functional $\vfi$. Then $(\ah,\deh)$ is again an algebraic quantum
group for which we can construct the dual $(\ahh,\dhh)$.
Theorem 4.12 of \cite{VD1} guarantees that $(\ahh,\dhh)$ is
isomorphic to $(A,\de)$:

\begin{theorem}
There exists an isomorphism of algebras $\Upsilon$ from
$A$ to $\ahh$ such that $\Upsilon(x)(\om) = \om(x)$ for every $x \in A$
and $\om \in \ah$. We have moreover that $\de \Upsilon = (\Upsilon \od
\Upsilon)\de$.
\end{theorem}

In the rest of this section, we use this theorem to identify $(A,\de)$
and $(\ahh,\dhh)$.

\begin{proposition}
Denote the universal corepresentation of $(A,\de)$ by $U$. Then
$\flip(U)$ is the universal corepresentation of $(\ah,\deh)$.
\end{proposition}
\begin{demo}
Choose $a,b \in A$. Then there exist $p_1,\ldots\!,p_n ,
q_1,\ldots\!,q_n \in A$ such that $b \ot a = \sum_{i=1}^n \de(p_i)(q_i
\ot 1)$. Looking at the proof of proposition \ref{prop5.1}, we see
that $U(b \ot a \vfi) = \sum_{i=1}^n q_i \ot p_i \vfi$. \newline This
implies that $\flip(U) \, (a \vfi \ot b) = \sum_{i=1}^n p_i \vfi \ot
q_i$. \ \ \ \ \ \ (*)

\medskip

Denote the universal corepresentation of $(\ah,\deh)$ by $V$. Using
the identification mentioned above, we get that $V$ is an element of
$M(\ah \ot A)$. Looking at section \ref{art5}, we consider $\ah \od A$
in this case as a subspace of $L(\ah)$. So we consider $V(a \vfi \ot
b)$ as an element of $L(\ah)$.

Therefore choose $\om \in \ah$.  Definition \ref{def5.1} implies that
$(V(a \vfi \ot b))(\om) = (\io \od b)(\deh(\om)(a \vfi \ot 1))$
\ \ (where $b$ is considered as an element of $\ah\,'$  ).

It is now easy to check for every $x \in A$ that
$$[(V(a \vfi \ot b))(\om)](x) = [(\io \od b)(\deh(\om)(a \vfi \ot
1))](x)= [\deh(\om)(a \vfi \ot 1)](x \ot b) \ . $$
Hence, remark \ref{rem2.1} implies that
\begin{eqnarray*}
& & [(V(a \vfi \ot b))(\om)](x) = (\om \od a \vfi)(\de(x)(b \ot 1))
= (\om \od \vfi)(\de(x) (b \ot a)) \\
& & \spat = \sum_{i=1}^n (\om \od \vfi)(\de(x p_i)(q_i \ot 1))
= \sum_{i=1}^n \vfi(x p_i) \om(q_i)
\end{eqnarray*}
where the left invariance was used once again in the last equality.
So
$$(V(a \vfi \ot b))(\om) = \sum_{i=1}^n \om(q_i) \, p_i \vfi
= [\sum_{i=1}^n p_i \vfi \ot q_i](\om) \ . $$
Therefore we get that $V(a \vfi \ot b) = \sum_{i=1}^n p_i \vfi \ot q_i$
which is equal to $\flip(U)\,(a \vfi \ot b)$ by equation (*).
\end{demo}

\medskip

\begin{corollary}
We have the equality $(\io \od \deh)(U) = U_{12} \, U_{13}$.
\end{corollary}
\begin{demo}
Because $\flip(U)$ is the universal corepresentation of $(\ah,\deh)$,
corollary \ref{cor5.1} implies that $(\deh \od \io)(\flip(U)) =
\flip(U)_{13} \, \flip(U)_{23}$. The corollary follows easily from
this equality.
\end{demo}

\end{document}